\documentclass[final]{siamltex}
\usepackage{amsmath,amsfonts,amssymb}
\usepackage{latexsym}
\usepackage{graphicx,overpic}
\usepackage{pgfplots}
\pgfplotsset{compat=newest}
\usepackage{color}
\usepackage{booktabs}
\usepackage{caption}
\usepackage{subcaption}
\usepackage{hyperref}
\usepackage{showlabels}
\usepackage{tikz}
\usepackage{lscape}
\usepackage{array}
\usepackage{colonequals}
\newcolumntype{C}[1]{>{\centering\let\newline\\\arraybackslash\hspace{0pt}}m{#1}}

\newcommand{\R}{\mathbb R}

\newcommand{\bA}{\mathbf A}

\newcommand{\bH}{\mathbf H}
\newcommand{\bI}{\mathbf I}

\newcommand{\bP}{\mathbf P}

\newcommand{\bN}{\mathbf N}

\newcommand{\bU}{\mathbf U}
\newcommand{\bV}{\mathbf V}

\newcommand{\ba}{\mathbf a}

\newcommand{\be}{\mathbf e}

\newcommand{\bk}{\mathbf k}
\newcommand{\bm}{\mathbf m}
\newcommand{\bn}{\mathbf n}

\newcommand{\bs}{\mathbf s}

\newcommand{\bu}{\mathbf u}
\newcommand{\bv}{\mathbf v}

\newcommand{\bx}{\mathbf x}

\newcommand{\bbf}{\mathbf f}
\newcommand{\blf}{\mathbf f}

\newcommand{\xx}{\boldsymbol{x}}

\newcommand{\T}{\mathcal T}

\newcommand{\divG}{{\mathop{\,\rm div}}_{\Gamma}}

\newcommand{\gradG}{\nabla_{\Gamma}}
\newcommand{\gradGh}{\nabla_{\Gamma_h}}
\newcommand{\curlG}{{\mathop{\,\rm curl}}_{\Gamma}}
\newcommand{\vcurlG}{{\mathop{\,\rm \mathbf{curl}}}_{\Gamma}}

\newcommand{\cT}{\mathcal T}

\newcommand{\nn}{\mathbb{N}}

\newcommand{\OGamma}{\Omega^\Gamma}

\renewcommand{\div}{\textrm{div}\ \!}

\DeclareGraphicsExtensions{.pdf,.eps,.ps,.eps.gz,.ps.gz,.eps.Y}

\DeclareMathOperator*{\RotG}{\mathbf{curl}_\Gamma \!}

\DeclareMathOperator*{\RotGh}{\mathbf{curl}_{\Gamma_h} \!}

\DeclareMathOperator*{\Div}{div \!}
\DeclareMathOperator*{\DivG}{div_\Gamma \!}
\DeclareMathOperator*{\DivGh}{div_{\Gamma_h} \!}

\DeclareMathOperator*{\nablaG}{\nabla_\Gamma \!}

\DeclareMathOperator*{\nablaGh}{\nabla_{\Gamma_h} \!}

\DeclareMathOperator*{\DeltaG}{\Delta_\Gamma \!}
\DeclareMathOperator*{\tr}{tr}

\newcommand{\vect}[1]{\boldsymbol{\mathbf{#1}}}

\newtheorem{assumption}{Assumption}[section]

\newtheorem{remark}{Remark}[section]
\begin{document}

\title{Finite element discretization methods for velocity-pressure and stream function formulations of surface Stokes equations}
\author{Philip Brandner\thanks{Institut f\"ur Geometrie und Praktische  Mathematik, RWTH-Aachen
University, D-52056 Aachen, Germany (brandner@igpm.rwth-aachen.de)} \and
Thomas Jankuhn\thanks{Institut f\"ur Geometrie und Praktische  Mathematik, RWTH-Aachen
University, D-52056 Aachen, Germany (jankuhn@igpm.rwth-aachen.de)} \and
Simon Praetorius\thanks{Institut f\"ur Wissenschaftliches Rechnen, TU Dresden, D-01062
Dresden, Germany (simon.praetorius@tu-dresden.de)} \and
Arnold Reusken\thanks{Institut f\"ur Geometrie und Praktische  Mathematik, RWTH-Aachen
University, D-52056 Aachen, Germany (reusken@igpm.rwth-aachen.de)} \and
Axel Voigt\thanks{Institut f\"ur Wissenschaftliches Rechnen, TU Dresden, D-01062 Dresden,
Germany, Center for Systems Biology Dresden (CSBD), Pfotenhauerstr. 108, 01307 Dresden, Germany and Cluster of Excellence - Physics of Life, TU Dresden, 01062 Dresden, Germany (axel.voigt@tu-dresden.de)}.
}
\maketitle

\begin{abstract}
In this paper we study parametric TraceFEM and parametric SurfaceFEM (SFEM) discretizations of a surface Stokes problem. These methods are applied both to the Stokes problem in velocity-pressure formulation and  in  stream function formulation. A class of  higher order methods is presented in a unified framework. Numerical efficiency aspects of the two formulations are discussed and a systematic comparison of TraceFEM and SFEM is given. A benchmark problem is introduced in which a scalar reference quantity is defined and numerically determined.
\end{abstract}


\section{Introduction}
Surface fluids arise in different applications such as emulsions, foams or biological membranes and can be modeled by surface (Navier-)Stokes equations (cf., e.g., \cite{scriven1960dynamics,slattery2007interfacial,arroyo2009,brenner2013interfacial,rangamani2013interaction,rahimi2013curved,Reuther_JCP_2016}). These equations constrain the velocity and pressure to a surface and, at least for stationary surfaces, enforce the velocity to be tangential to the surface, which leads to a tide coupling with geometric properties of the surface and new physical phenomena. Despite the apparent practical relevance, there has been only recently a strongly growing mathematical interest in modeling of surface fluids, e.g., \cite{arroyo2009,nitschke2012finite,Jankuhn1,Kobaetal_QAM_2017,Koba_2017,miura2017singular,Nitschkeetal_arXiv_2018,ReVo2015,ReVo2018,ToMiAr2019,reuther2020numerical} and their numerical simulation, e.g., \cite{nitschke2012finite,barrett2014stable,ReVo2015,Nitschke_DEC_2017,Reusken_Streamfunction,reuther2018solving,fries2018higher,olshanskii2018finite,BrandnerReusken2019,olshanskii2019penalty,Bonito2019a,OlshanskiiZhiliakov2019,ToMiAr2019,Lederer2019,reuther2020numerical,ToHu2021}. Surface (Navier-)Stokes equations are also studied as an interesting mathematical problem on its own, e.g.,  \cite{ebin1970groups,Temam88,taylor1992analysis,arnol2013mathematical,mitrea2001navier,arnaudon2012lagrangian}.

In the discretization of surface (Navier-)Stokes equations several issues occur, which are not present for the (Navier-)Stokes equations in the standard Euclidean space. For example, there are difficulties related to the approximation of the surface $\Gamma$ and of several  quantities associated with the geometry such as  covariant derivatives and  curvature terms. Another difficulty is to ensure tangency of the velocity field. Most of the cited approaches enforce the tangential condition weakly: using a Lagrange multiplier (cf. \cite{fries2018higher,jankuhn2019}) or a penalty term (cf. \cite{reuther2018solving,olshanskii2019penalty}). Such approaches are applied both in trace finite element methods (TraceFEM) and in  surface finite element methods (SFEM). In \cite{Lederer2019,Bonito2019a} an alternative SFEM is considered, in which a Piola transformation for the construction of divergence-free tangential finite elements is introduced. In this paper we restrict to the most popular technique for handling the tangential condition, namely the penalty method. Instead of treating the (Navier-)Stokes equations in the velocity and pressure variables one can also use a stream function formulation (cf. \cite{nitschke2012finite,ReVo2015,ReVo2018,Reusken_Streamfunction,Grossetal_JCP_2018,BrandnerReusken2019,ToMiAr2019}). The approach has the advantage that only scalar quantities have to be considered. The velocity can be approximated from the computed stream function. In this setting there is no difficulty concerning tangency of the velocity field.

In this paper we compare two discretization methods for the surface Stokes equations, namely the parametric TraceFEM and the parametric SFEM. We consider both a formulation in the velocity and pressure variables and a stream function formulation. For TraceFEM the first formulation is treated in \cite{jankuhn2020error} and the second one is based on \cite{BrandnerReusken2019}.  For SFEM the first formulation extends the approach in \cite{reuther2018solving,fries2018higher} and the second formulation is based on \cite{nitschke2012finite}. We outline the key components of these methods and discuss further related literature. For the formulation in velocity and pressure variables we use generalized Taylor-Hood elements $\textbf{P}_k-P_{k-1}$, $k\geq 2$, defined on the bulk mesh (TraceFEM) or the surface mesh (SFEM). A consistent penalty approach is used for both methods to satisfy the tangential constraint weakly. Higher order TraceFEM is obtained using the parametric finite element approach introduced for scalar problems in \cite{lehrenfeld2016high}.
For this TraceFEM  a  stability and discretization error analysis including geometrical errors is presented in \cite{jankuhn2020error}. A $\textbf{P}_1-P_1$ variant of the SFEM was first introduced in \cite{reuther2018solving} and numerical simulation results with the $\textbf{P}_k-P_{k-1}$, $k\geq 2$, Taylor-Hood pairs are given in \cite{fries2018higher}. The  higher order parametric SFEM that we present extends the approaches used in \cite{DEreview,Nestleretal_arxiv_2018,demlow2009higher}. Error analysis of the SFEM approach for surface Stokes problems are not available in the literature. An error analysis of this method for a surface vector-Laplace equation is presented in \cite{hansbo2019analysis}. Related to the stream function formulation we note the following. This approach requires the surface to be simply connected. In the fields of applications mentioned above, one often deals with smooth simply connected surfaces without boundary. In such a setting there usually are no difficulties related to regularity or boundary conditions and the stream function formulation may be an attractive alternative to the formulation in velocity-pressure variables, as already indicated in \cite{nitschke2012finite}. In \cite{Reusken_Streamfunction} fundamental properties of the surface stream function formulation, e.g. with respect to well-posedness and  relations to a surface Helmholtz decomposition, are derived. In both papers \cite{nitschke2012finite,Reusken_Streamfunction} the resulting fourth order scalar surface partial differential equation for the stream function is reformulated as a coupled system of two second order equations, which is a straightforward generalization to surfaces of the classical Ciarlet-Raviart method \cite{Ciarlet1974} in  Euclidean space. As the equations are scalar-valued they can be discretized by established finite element methods for scalar-valued surface partial differential equations, such as TraceFEM \cite{olshanskii2016trace}, SFEM \cite{DEreview} or diffuse interface approximations \cite{Raetz2006}; cf. also the overview paper \cite{Bonito2019}. In \cite{BrandnerReusken2019} an error analysis of the TraceFEM for the stream function formulation of the surface Stokes equations is presented. The main new contributions of the paper are the following:
\begin{itemize}
 \item we present a \emph{general methodology for optimal higher order TraceFEM and SFEM}. Several key ingredients are known from the literature, e.g.  the higher order surface approximation methods introduced in \cite{demlow2009higher,lehrenfeld2016high}. These are combined with suitable parametric finite element spaces and  methods for computing ``sufficiently accurate'' normal and Gauss curvature approximations.
 \item we present a \emph{systematic comparison of the velocity-pressure and the stream function formulations} of surface Stokes. Both approaches are natural ones, but so far they have not been compared for surface Stokes equations.
 \item we present a \emph{systematic comparison of TraceFEM and SFEM}. We compare specific measures of complexity of the two methods and present numerical simulation results that allow comparison of the two methods.
 \item we introduce a \emph{benchmark problem for surface Stokes} equations. We define a scalar quantity (related to the distance between a vortex in the solution and a maximal curvature location on the surface) that is determined with our simulation codes. Since we have different formulations and different finite element methods that are implemented in different codes, we can determine with high reliability the accuracy of the computed reference quantity.
\end{itemize}

The remainder of the paper is organized as follows. In Section~\ref{sec:continuous_problem} we introduce surface differential operators and recall the well-posed weak formulations for the surface Stokes equations of both formulations. Parametric approximations for surfaces are explained in Section \ref{sectparametric} and TraceFEM and SFEM approaches for both formulations are presented in Section \ref{sec:discretization_methods}. Finally, in Section~\ref{sec:numerical_experiments} we compare the two problem formulations and the two discretization methods numerically. The section also contains the benchmark problem.

\section{Continuous problem} \label{sec:continuous_problem}
We consider a smooth hypersurface $\Gamma$ without boundary and a polygonal domain $\Omega \subset \mathbb{R}^3$ with $\Gamma \subset \Omega$. Let $d$ denote the signed distance function to $\Gamma$ which is negative in the interior of $\Gamma$. For $\delta>0$ we define the neighborhood $U_\delta  \colonequals  \left\lbrace \xx \in \mathbb{R}^3 \mid \vert d(\xx) \vert < \delta \right\rbrace$ of $\Gamma$. For $\delta > 0$ sufficiently small and $\xx \in U_\delta$ we define $\bn(\xx) = \nabla d(\xx)$ (for $\xx \in \Gamma$ this is the outward pointing unit normal), the orthogonal projection  $\bP = \bP(\xx) \colonequals  \bI - \bn(\xx)\bn(\xx)^T$, the closest point projection $\pi(\xx) = \xx - d(\xx)\bn(\xx)$ and the Weingarten map $\bH(\xx) = \nabla^2d(\xx)$. We assume that $\delta$ is sufficiently small such that the decomposition $\xx = \pi(\xx) + d(\xx) \bn(\xx)$ is unique for all $\xx \in U_{\delta}$. Let $\psi^e(\xx) \colonequals \psi(\pi(\xx))$ and $\bv^e(\xx) \colonequals \bv(\pi(\xx))$ for $\xx \in U_{\delta}$ be the constant normal extension for scalar functions $\psi \colon \Gamma \to \mathbb{R}$ and vector functions $\bv \colon \Gamma \to \mathbb{R}^3$, respectively. The tangential surface derivatives for scalar functions $\psi \colon \Gamma \to \mathbb{R}$ and vector functions $\bv \colon \Gamma \to \mathbb{R}^3$ are defined by
\begin{equation}\label{eq:surface_gradients} \begin{split}
  \nabla_{\Gamma} g(\xx) &= \bP(\xx)\nabla g^e(\xx), \quad \xx \in \Gamma, \\
  \nabla_{\Gamma} \bv(\xx) &= \bP(\xx)\nabla \bv^e(\xx) \bP(\xx), \quad \xx \in \Gamma.
\end{split}
\end{equation}
To simplify the notation we often drop the argument $\xx$. For a vector field $\bu$ the (infinitesimal) deformation tensor is given by
\begin{equation*}
  E_s(\bu) \colonequals \frac{1}{2} \left( \gradG \bu + \gradG \bu^T \right).
\end{equation*}
Let $\be_i$ be the $i$th basis vector in $\mathbb{R}^3$. We define the surface divergence operator for vector-valued functions $\bu \colon \Gamma \to \mathbb{R}^3$ and tensor-valued functions $\bA \colon \Gamma \to \mathbb{R}^{3\times3}$ by
\begin{equation*} \begin{split}
  \divG \bu & \colonequals  \textrm{tr} (\gradG \bu),  \\
  \divG \bA & \colonequals  \left( \divG (\be_1^T\bA), \divG (\be_2^T\bA),\divG (\be_3^T\bA)  \right)^T.
\end{split}
\end{equation*}
Note that in the literature there are other definitions of the surface divergence, in which an additional surface-projection is included, cf. \cite{reuther2018solving}.  The surface curl operators are defined by
\begin{align*}
  \curlG \bu \colonequals \DivG(\bu \times \bn), \quad \bu \in C^1(\Gamma)^3, \\
  \vcurlG \phi \colonequals  \bn \times \gradG \phi, \quad \phi \in C^1(\Gamma).
\end{align*}

For a given force vector $\bbf \in L^2(\Gamma)^3$ with $\bbf \cdot \bn =0$ we consider the following \emph{surface Stokes problem}: Determine $\bu \colon \Gamma \to \R^3$ with $\bu \cdot \bn = 0$ and $p \colon \Gamma \to \R$ with $\int_\Gamma p\, ds=0$ such that
\begin{equation} \label{eq:strong} \begin{split}
- \bP \divG (E_s(\bu)) + \bu + \gradG p &= \bbf \qquad \text{on } \Gamma,  \\
\divG \bu &= 0 \qquad \text{on } \Gamma.
\end{split}
\end{equation}
The zero order term on the left-hand side is added to avoid technical details related to the kernel of the tensor $E_s$, also called the space of Killing vector fields. Below we recall two variational formulations of the surface Stokes problem \eqref{eq:strong}.

\begin{remark} \label{rem:strong_alternative} Alternative formulations for the surface Stokes problem \eqref{eq:strong} are obtained by using the identities $2 \bP \divG (E_s(\bu)) = \boldsymbol{\Delta}_{\Gamma}^B \bu + K \bu = - \boldsymbol{\Delta}_{\Gamma}^{dR} \bu + 2 K \bu$, where $K$ denotes the Gaussian curvature of the surface $\Gamma$, $\boldsymbol{\Delta}_\Gamma^B \bu = \bP \divG \nabla_{\Gamma} \bu$ is the Bochner Laplacian and  $\boldsymbol{\Delta}_{\Gamma}^{dR} \bu = - (\vcurlG \curlG + \nabla_{\Gamma} \divG) \bu$ the Laplace-deRham operator, see \cite{AMR88}.
\end{remark}

\subsection{Variational formulation in $\bu$-$p$ variables}
We recall a standard weak formulation of the surface Stokes problem in velocity--pressure variables. For this we need the  surface Sobolev space of weakly differentiable vector-valued functions, denoted by $H^1(\Gamma)^3$, with norm $\Vert \bu \Vert_{H^1(\Gamma)}^2  \colonequals  \int_{\Gamma} \Vert \bu(s) \Vert_2^2 + \Vert \nabla \bu^e(s) \Vert_2^2 \, ds$. The corresponding subspace of \emph{tangential} vector fields is denoted by
$
  \bH_t^1(\Gamma) \colonequals \{ \, \bu \in H^1(\Gamma)^3~|~  \bu \cdot \bn=0 \quad \text{a.e. on}~\Gamma\,\}\,.
$
A vector $\bu \in H^1(\Gamma)^3$ can be orthogonally decomposed into a tangential and a normal part. We use the notation:
\begin{equation*}
  \bu = \bP \bu + (\bu \cdot \bn)\bn = \bu_T + u_N\bn\,.
\end{equation*}
For $\bu, \bv\in H^1(\Gamma)^3$ and $p \in L^2(\Gamma$) we introduce the bilinear forms
\begin{align}
  \ba(\bu, \bv) &\colonequals \int_\Gamma E_s(\bu) : E_s(\bv) \, ds + \int_\Gamma \bu\cdot\bv \, ds\,, \label{blfa} \\
  b_T(\textbf{u}, p) &\colonequals -\int_\Gamma p \divG \bu_T \, ds\,. \label{blfb}
\end{align}
Note that in the definition of $b_T(\bu,p)$ only the \emph{tangential} component of $\bu$ is used, i.e., $b_T(\bu,p) = b_T(\bu_T,p)$ for all $\bu \in H^1(\Gamma)^3$, $p\in L^2(\Gamma)$. This property motivates the notation $b_T(\cdot,\cdot)$ instead of $b(\cdot,\cdot)$. If $p$ is from $H^1(\Gamma)$, then integration by parts yields
\begin{equation}\label{Bform}
  b_T(\bu,p) = \int_\Gamma \bu_T\cdot \gradG p \, ds = \int_\Gamma \bu \cdot \gradG p \, ds\,.
\end{equation}
We introduce the following variational formulation: Determine $(\bu_T, p)  \in \bH_t^1(\Gamma) \times L^2_0(\Gamma)$ such that
\begin{equation} \label{contform} \begin{aligned}
  \ba(\bu_T, \bv_T) + b_T(\bv_T, p) &= (\bbf, \bv_T)_{L^2(\Gamma)} && \text{for all } \bv_T \in \bH_t^1(\Gamma), \\
b_T(\bu_T,q) &= 0 && \text{for all } q \in L^2(\Gamma)\,.
\end{aligned}
\end{equation}
This is a \emph{well-posed variational formulation} of the surface Stokes problem \eqref{eq:strong}, cf.~ \cite{Jankuhn1}. The unique solution is denoted by $(\bu_T^*,p^*)$. For the discretization, we need the bilinear form $\ba_T(\bu,\bv) \colonequals \ba(\bP \bu, \bP \bv) = \ba(\bu_T,\bv_T)$. Using the identity $E_s(\bu )= E_s(\bu_T) + u_N \bH$ we get
\begin{equation} \label{idneeded}
  \ba_T(\bu,\bv) = \int_\Gamma \big( E_s(\bu)- u_N \bH \big):\big(E_s(\bv)- v_N \bH \big)  + \bu_T \cdot \bv_T \, ds\,.
\end{equation}

\subsection{Variational formulation in stream function variable}
We recall the stream function formulation of the surface Stokes problem \cite{Reusken_Streamfunction}. For its derivation we need the following assumption:

\begin{assumption} \label{ass} In the remainder we assume that $\Gamma$ is simply connected and sufficiently smooth, at least $C^3$.
\end{assumption}

We introduce the spaces $
  \bH_{t,\div}^1 \colonequals \{\, \bu \in \bH_t^1(\Gamma)~|~\DivG \bu =0 \,\}$, $
  H_\ast^k(\Gamma) \colonequals \{\, \psi \in H^k(\Gamma)~|~\int_\Gamma \psi \, ds =0\,\}\,,
$
and the bilinear form
\begin{equation*}
  a(\phi,\psi) \colonequals \int_\Gamma \tfrac12 \DeltaG \phi \DeltaG \psi + (1- K)\gradG \phi \cdot \gradG \psi \, ds,
\end{equation*}
with $K$  the Gaussian curvature of the surface $\Gamma$. The following result is derived in \cite{Reusken_Streamfunction} for the case without a zero order term in \eqref{eq:strong}. Without any significant changes, this derivation also applies to \eqref{eq:strong}.

\begin{theorem} \label{thm:stream}
  Let $\bu_T^\ast \in \bH_{t,\Div}^1$ be the unique solution of \eqref{contform} and $\psi^\ast \in H_\ast^1(\Gamma)$ its unique stream function, i.e., $\bu^\ast= \vcurlG \psi^\ast$. This $\psi^\ast$ is the unique solution of the following stream function problem: Determine $\psi \in H_\ast^2(\Gamma)$ such that
  \begin{equation} \label{eq:weakstream}
    a(\psi,\phi)= (\blf, \vcurlG \phi)_{L^2(\Gamma)} \quad \text{for all }\phi \in H_\ast^2(\Gamma)\,.
  \end{equation}
\end{theorem}

In view of a finite element discretization it is convenient to reformulate the fourth order surface equation \eqref{eq:weakstream} as a coupled system of two second order problems, introducing the vorticity $\phi$, similar as for the classical two-dimensional Stokes problem. We define the bilinear forms:
\[
  m(\xi,\eta) :=  \int_{\Gamma} \xi \eta \, ds\,,~~
  \ell(\xi,\eta):= \int_{\Gamma} \nablaG \xi \cdot  \nablaG \eta \, ds\,, ~~
  \ell_K(\xi,\eta) := 2 \int_{\Gamma} (1-K) \nablaG \xi \cdot  \nablaG \eta \, ds\,,
\]
and the linear functional
$
  g(\xi) :=  -2 \int_{\Gamma} \blf \cdot \RotG \xi \, ds\,.
$
The coupled second order system is as follows: Determine $\psi \in H^1_*(\Gamma)$, $\phi \in H^1(\Gamma)$ such that
\begin{equation}\label{NotesAR_8}
\begin{aligned}
  m(\phi,\eta) + \ell(\psi, \eta) &= 0      &&\text{for all } \eta \in H^1(\Gamma)\,,\\
  \ell(\phi,\xi) - \ell_K(\psi, \xi) &= g(\xi) &&\text{for all } \xi \in H^1(\Gamma)\,.
\end{aligned}
\end{equation}
In \cite{Reusken_Streamfunction} it is shown that this problem has a unique solution  $\psi = \psi^\ast$, $\phi = \phi^\ast = \DeltaG \psi^\ast$, with $\psi^\ast$ being the unique solution of \eqref{eq:weakstream}. In the remainder we denote by $\psi^\ast$ and $\phi^\ast$ the unique solution of \eqref{NotesAR_8}. Below we introduce  finite element discretization methods that  are based on \eqref{NotesAR_8}.

We briefly address natural variational formulations that can be used to determine the velocity solution $\bu_T^*$ and the pressure solution $p^*$, given the stream function solution $\psi^\ast$. Based on the relation $\bu_T^\ast =\RotG \psi^\ast = \bn \times \nablaG \psi^\ast$ we introduce a well-posed variational formulation for the velocity reconstruction: Determine $\bu \in L^2(\Gamma)^3$ such that
\begin{align}\label{Variationsformulierung_u}
  \int_{\Gamma} \bu \cdot \bv \, ds = \int_{\Gamma} \left( \bn \times \nablaG \psi^\ast \right) \cdot \bv \, ds  \quad \text{for all } \bv \in L^2(\Gamma)^3\,.
\end{align}
The unique solution of this problem coincides with $\bu_T^\ast$. For the pressure reconstruction we introduce the variational problem: Determine $p\in H^1_\ast(\Gamma)$ such that
\begin{equation} \label{Variationsformulierung_p}
  \int_\Gamma \nabla_\Gamma p \cdot \nabla_\Gamma \xi \, ds =  \int_\Gamma \left(  K \RotG \psi^\ast + \blf  \right) \cdot \nabla_\Gamma \xi \, ds \quad \text{for all } \xi \in H^1(\Gamma)\,.
\end{equation}
In \cite{BrandnerReusken2019} it is shown that the pressure solution $p^\ast$ coincides with the unique solution of the Laplace-Beltrami problem \eqref{Variationsformulierung_p}. The variational problems \eqref{Variationsformulierung_u} and \eqref{Variationsformulierung_p} can be used in finite element reconstruction methods for the velocity and pressure solutions, respectively.

\section{Parametric finite elements for  surface approximation} \label{sectparametric}
For a higher-order finite element discretizations of the surface Stokes problem we need a more accurate than piecewise linear surface approximation. Different techniques for constructing higher order surface approximations are available in the literature, e.g., \cite{demlow2009higher,Bonito2019,lehrenfeld2016high,grande2017higher,FriesOmerovicSchoellhammerStanford2016,FriesSchoellhammer2017,FriesOmerovic2016,praetorius2020dune}, for both, TraceFEM and SFEM. In the subsections \ref{sectS1} and \ref{sectS2} below we briefly recall two known techniques.

\subsection{Surface approximation for the TraceFEM} \label{sectS1}
We outline the technique introduced in \cite{lehrenfeld2016high}, for which it is essential that $\Gamma$ is characterized as the zero level of a smooth level set function $\varphi \colon U_\delta \to \mathbb{R}$, i.e., $\Gamma = \{ \xx \in \Omega \mid \varphi(\xx) = 0 \}$. We do not assume the level-set function to be close to a distance function but to have the usual properties of a level-set function: $\Vert \nabla \varphi(\xx) \Vert \sim 1$, $\Vert \nabla^2 \varphi(\xx) \Vert \leq C$ for all $\xx \in U_\delta$. Let $\{ \mathcal{T}_h \}_{h>0}$ be a family of shape regular tetrahedral triangulations of $\Omega$. By $V_h^k$ we denote the standard finite element space of continuous piecewise polynomials of degree $k$. The nodal interpolation operator in $V_h^k$ is denoted by $I^k$. As input for a parametric mapping we need an approximation of $\varphi$. The construction of the geometry approximation will be based on a level set function approximation  $\varphi_h \in V_h^{k}$. We assume that for this approximation the error estimate
\begin{equation} \label{phibound}
  \max_{T \in \mathcal{T}_h} \vert \varphi_h - \varphi \vert_{W^{l,\infty}(T\cap U_\delta)} \leq C\, h^{k+1-l}, \quad 0\leq l \leq k+1\,,
\end{equation}
is satisfied. Here, $\vert \cdot \vert_{W^{l,\infty}(T\cap U_\delta)}$ denotes the usual semi-norm in the Sobolev space $W^{l,\infty}(T\cap U_\delta)$ and the constant $c$ depends on $\varphi$ but is independent of $h$. The zero-level set of the finite element function $\varphi_h$ \emph{implicitly} defines an approximation of the interface, on which, however, numerical integration is hard to realize for $k \geq 2$.  With the piecewise \emph{linear} nodal interpolation of $\varphi_h$, which is denoted by $\hat{\varphi}_h = I^1\varphi_h$, we define the low order geometry approximation
$
  \Gamma^{\text{lin}}  \colonequals  \{ \xx \in \Omega \mid \hat{\varphi}_h (\xx) = 0\}\,$. This piecewise planar surface approximation in general is very shape \emph{ir}regular.
The tetrahedra $T \in \mathcal{T}_h$ that have a nonzero intersection with $\Gamma^{\text{lin}}$ are collected in the set denoted by $\T_h^\Gamma$. The domain formed by all tetrahedra in $\T_h^\Gamma$ is denoted by $\OGamma_h \colonequals  \{ x \in T \mid T \in \T_h^\Gamma \}$. Let $\Theta_h \in  \big({V_h^{k}}_{|\OGamma_h}\big)^3$ be the mesh transformation of order $k$ as defined in \cite{lehrenfeld2016high,grande2017higher}. An approximation of $\Gamma$ is defined by
\begin{equation}\label{eq:discrete-surface}
  \Gamma_{h} \colonequals \Theta_h(\Gamma^{\text{lin}}) = \left\lbrace \xx \mid \hat{\varphi}_h(\Theta_{h}^{-1}(\xx)) = 0 \right\rbrace.
\end{equation}
In \cite{lehrenfeld2017analysis} it is shown that (under certain reasonable smoothness assumptions) the estimate
$
  {\rm dist}(\Gamma_{h}, \Gamma) \lesssim h^{k+1}
$
holds. Hence, the parametric mapping $\Theta_{h}$ indeed yields a higher order surface approximation. Here and further in the paper we write $x \lesssim y$ to state that there exists a constant $C>0$, which is independent of the mesh parameter $h$ and the position of $\Gamma$ in the background mesh, such that the inequality $x \leq C\,y$ holds. We denote the transformed cut mesh domain by $\Omega^{\Gamma}_\Theta  \colonequals  \Theta_h(\Omega^\Gamma_h)$ and apply to $V_h^{k}$ the transformation $\Theta_h$ resulting in the isoparametric spaces (defined on $\Omega^{\Gamma}_\Theta$)
\begin{equation} \label{FEspace1}
  V_{h,\Theta}^{k} \colonequals \left\lbrace v_h \circ \Theta_h^{-1} \mid v_h \in {V_h^{k}}|_{\Omega^{\Gamma}_h} \right\rbrace, \quad
  \bV_{h,\Theta}^{k} \colonequals (V_{h,\Theta}^{k})^3\,.
\end{equation}

The following lemma, taken from \cite{grande2017higher}, gives an approximation error for the easy to compute normal approximation $\bn_h$, which is used in the methods introduced below.

\begin{lemma} \label{lemmanormals}
For $\xx \in T \in \mathcal{T}^\Gamma_h$ define
\begin{equation*}
  \bn_{\textrm{lin}} = \bn_{\textrm{lin}}(T) \colonequals \frac{\nabla \hat{\varphi}_h(\xx)}{\Vert \nabla \hat{\varphi}_h(\xx)\Vert_2} = \frac{\nabla \hat{\varphi}_{h|T}}{\Vert \nabla \hat{\varphi}_{h|T}\Vert_2}, \quad
  \bn_h (\Theta_h(\xx)) \colonequals \frac{D\Theta_h(\xx)^{-T} \bn_{\textrm{lin}}}{\Vert D\Theta_h(\xx)^{-T} \bn_{\textrm{lin}} \Vert_2}\,.
\end{equation*}
Restricted to the surface the approximations $\bn_{\textrm{lin}}$ and $\bn_h $ are normals on $\Gamma_{\textrm{lin}}$ and $\Gamma_h$, respectively. Furthermore, the following holds:
\begin{equation*}
  \Vert \bn_{h} - \bn \Vert_{L^\infty(\Omega^{\Gamma}_\Theta)} \lesssim h^{k}.
\end{equation*}
\end{lemma}

\subsection{Surface approximation for the SFEM}\label{sectS2}
The method that we outline in this section is based on \cite{demlow2009higher}. We assume that  $\Gamma$ is given as the zero level of a level set function. This assumption is not essential. For constructing a higher order surface approximation we start from a piecewise planar surface approximation $\Gamma^{\text{lin}}$. However, different from the approach discussed in section~\ref{sectS1}, it is essential that  this initial approximation is  shape-\emph{regular}. In our implementation we use  quasi-uniform triangulations $\Gamma^{\text{lin}}$.
  Such an approximate surface triangulation can be generated by an approach based on optimizing the local element quality using vertex-motion and edge flipping, see, e.g., \cite{persson2004distmesh,persson2004thesis}, or using a mesh coarsening approach that optimizes the edge lengths and angles during mesh simplification and decimation, see, e.g., \cite{valette2004acvd,meshconv2020}.

The construction of higher-order surface triangulations then follows the approach described in \cite{demlow2009higher} which is implemented in \cite{praetorius2020dune} for the \textsc{Dune} discretization framework. Let $\mathcal{S}^{\text{lin}}_h$ be the shape-regular (surface) triangulation of $\Gamma^{\text{lin}} = \bigcup_{\hat{S}\in\mathcal{S}^{\text{lin}}_h} \hat{S}$ with each element $\hat{S}$  parametrized over a reference domain $\Lambda\subset\mathbb{R}^2$ by $F_{\hat{S}}(\Lambda)=\hat{S}$.
A surface-mesh transformation is based on a piecewise Lagrange polynomial interpolation of the closest point projection $\pi:U_\delta\to\Gamma$ using local Lagrange basis functions $\{\vartheta^i\}_{1\ldots n_k}$ of order $k$ defined on the reference domain $\Lambda$. Let such an interpolation function be denoted by $\pi_h\big|_{\hat{S}} = I^k (\pi\circ F_{\hat{S}})$ with $I^k$ the $k$-th order Lagrange interpolation operator on $\Lambda$, i.e.,
\begin{equation*}
  \bar{\pi}_{\hat{S},h}(\lambda)\colonequals \pi_h(F_{\hat{S}}(\lambda)) = \sum_{i=1}^{n_k} \pi(F_{\hat{S}}(\lambda_i)) \vartheta^i(\lambda)\quad\text{ for }\lambda\in \Lambda\,,
\end{equation*}
with $\lambda_i\in\Lambda$ the local Lagrange nodes on $\Lambda$ corresponding to the local Lagrange basis function $\vartheta^i$. Mapping the nodes of all elements of $\mathcal{S}^{\text{lin}}_h$ yields the piecewise polynomial surface of order $k$,
\begin{equation*}
  \Gamma_h
    \colonequals \pi_h(\Gamma^\text{lin})
    = \bigcup_{\hat{S}\in\mathcal{S}^{\text{lin}}_h} \{\pi_h(\hat{\xx})\;\vert\;\hat{\xx}\in \hat{S}\}
    = \bigcup_{\hat{S}\in\mathcal{S}^{\text{lin}}_h} \pi_h(\hat{S})
    \equalscolon \bigcup_{S\in\mathcal{S}_h} S\,.
\end{equation*}
Note that this construction of the discrete surface $\Gamma_h$ is different from \eqref{eq:discrete-surface} but also leads to a piecewise polynomial approximation of $\Gamma$. In the discussion of the methods and numerical results below, the corresponding form of the discrete surface has to be taken into account.

The Lagrange finite element space of order $m$ on the piecewise flat triangulated surface $\Gamma^{\text{lin}}$, defined by $\hat{V}_h^m(\Gamma^{\text{lin}}) = \{\hat{v}\in C^0(\Gamma^{\text{lin}})\;\vert\;\hat{v}|_{\hat{S}}\in\mathbb{P}_m\;\forall {\hat{S}}\in\mathcal{S}^{\text{lin}}_h\}$, induces  a corresponding Lagrange finite element space on the polynomial surface $\Gamma_h$, by lifting the functions to the curved elements. In this paper we restrict to the isoparametric case $m=k$ and thus we get the spaces
\begin{equation*}
  \tilde V_{h,\pi_h}^k(\Gamma_h) \colonequals \big\{\hat{v}\circ \pi_h^{-1}\;\vert\;\hat{v}\in\hat{V}_h^k \big\}\,,~~\tilde \bV_{h,\pi_h}^k\colonequals \big(\tilde V_{h,\pi_h}^k\big)^3,
\end{equation*}
which can be compared to the spaces defined in \eqref{FEspace1}.

In the SFEM for the surface Stokes equations, surface normals and the Weingarten map for the parametric surface $\Gamma_h$ are required. These can be obtained from the derivatives of the polynomial projection function $\pi_h$, see also \cite{praetorius2020dune}. We denote by $\bar{J}^i\colonequals\nabla_{\lambda}\vartheta^i$ the local basis function Jacobians. The Jacobian of $\pi_h$ and the normal vectors on $\Gamma_h$ are then given by
\begin{align*}
  \bar{J}_{S}(\lambda) &\colonequals \sum_{i=1}^{n_k} \pi(F_{\hat{S}}(\lambda_i))\otimes \bar{J}^i(\lambda),\quad\bar{\bn}_{S}(\lambda) = \frac{\bar{\bN}_{S}(\lambda)}{\|\bar{\bN}_{S}(\lambda)\|}\,\text{ for }\lambda\in\Lambda
\end{align*}
 with $\bar{\bN}_{S}(\lambda) \colonequals \bar{J}_{S}(\lambda)_{\cdot,1} \times \bar{J}_{S}(\lambda)_{\cdot,2}$ the cross product of the columns of $\bar{J}_{S}(\lambda)$. We identify $\bn_h(\xx) = \bn_h(\bar{\pi}_{\hat{S},h}(\lambda))\equiv \bar{\bn}_{S}(\lambda)$ for $\xx=\bar{\pi}_{\hat{S},h}(\lambda)\in S$ and $\bN_h\circ\bar{\pi}_{\hat{S},h} \equiv \bar{\bN}_{S}$, analogously. The approximate Weingarten map $\bH_h(\xx)=\nabla_{S}\bn_h(\xx)$ for $\xx\in S$ then follows by chain rule using the surface derivatives \eqref{eq:surface_gradients} of $\bN_h$:
\begin{equation*}
  \bH_h \colonequals \frac{\nabla_{S}\bN_h}{\|\bN_h\|}\,,
\end{equation*}
locally, inside each element $S\in\mathcal{S}_h$. In \cite{demlow2009higher} the following estimates for errors in the normal vector and Weingarten map are proven:
\begin{equation} \label{estDemlow}
  \|\bn_h - \bn\circ\pi\|_{L^\infty(\Gamma_h)} \lesssim h^k\,,\quad
  \|\bH_h - \bH\circ\pi\|_{L^\infty(\Gamma_h)} \lesssim h^{k-1}\,.
\end{equation}

All these surface approximations are based on an initial linear approximation $\Gamma^{\text{lin}}$ of $\Gamma$ and an interpolation of the exact closest point projection $\pi$, which is not always available directly. Thus, it needs to be computed numerically. In \cite{persson2004thesis,demlow2007adaptive,nitschke2014thesis} some iterative schemes are discussed to evaluate the closest point projection in the  neighborhood of $\Gamma$. We follow the iterative approach introduced in \cite{demlow2007adaptive} in our numerical experiments, see also \cite{praetorius2020dune}.

\section{Discretization methods for the surface Stokes equations} \label{sec:discretization_methods}
In this section we outline the two discretization methods TraceFEM and SFEM for the discretization of the surface Stokes variational problems \eqref{contform} and \eqref{NotesAR_8}, and the variational formulations for the reconstruction of the velocity and pressure \eqref{Variationsformulierung_u}, \eqref{Variationsformulierung_p}.

\subsection{Trace finite element method} \label{sec:trace_fem}
Since the TraceFEM is a geometrically \emph{un}fitted finite element method, we need a stabilization that eliminates instabilities caused by the small cuts. We use the so-called ``normal derivative volume stabilization'' \cite{burmanembedded, grande2017higher}:
\begin{equation*}
  \bs_h(\bu,\bv) \colonequals  \rho_{\bu} \int_{\Omega_{\Theta}^{\Gamma}} (\nabla \bu \bn_h) \cdot (\nabla \bv \bn_h)  \, dx\,, \quad
  s_h(p,q)  \colonequals  \rho_p \int_{\Omega_{\Theta}^{\Gamma}} (\bn_h\cdot \nabla p) (\bn_h \cdot \nabla q) \, dx\,,
\end{equation*}
with parameters $\rho_{\bu}$ and $\rho_p$ specified below, cf. Table~\ref{tabpar}.

\subsubsection{Discretization of variational formulation in $\bu$-$p$ variables}\label{sec:trace_fem_th}
Based on the parametric finite element spaces $\bV_{h,\Theta}^{k}$ and $V_{h,\Theta}^{k}$ we introduce  the $\vect P_k$-$P_{k-1}$ pair of parametric trace Taylor-Hood elements:
\begin{equation} \label{TaylorHood}
  \bU_h  \colonequals  \bV_{h,\Theta}^{k}, \qquad Q_h  \colonequals  V_{h,\Theta}^{k-1} \cap L^2_0(\Gamma_h), \quad k \geq 2\,.
\end{equation}
Note that  the polynomial degrees, $k$ and $k-1$, for the velocity and pressure approximation are different, but both spaces $\bU_h$ and $Q_h$ use the same parametric mapping based on polynomials of degree $k$. Since the pressure approximation uses $H^1$ finite element functions we can use the integration by parts \eqref{Bform} (with $\Gamma$ replaced by $\Gamma_h$).
We introduce discrete variants of the bilinear forms $\ba_T(\cdot,\cdot)$, cf. \eqref{idneeded},  and $b_T(\cdot,\cdot)$. We define, with $\bP_h = \bP_h(\xx) \colonequals \bI- \bn_h(\xx)\bn_h(\xx)^T$, $\xx \in \Omega_{\Theta}^{\Gamma}$, $\bu,\bv \in \bU_h$, $q \in Q_h$:
\begin{align*}
  \gradGh q   &\colonequals \bP_h \nabla q\,, \quad
  \gradGh \bu  \colonequals \bP_h \nabla \bu \bP_h\,, \\
  E_h(\bu)    &\colonequals \tfrac12 \big(\gradGh \bu + \gradGh \bu^T \big)\,, \quad
  E_{T,h}(\bu) \colonequals E_h(\bu) - (\bu \cdot \bn_h) \bH_h\,, \\
  \ba_{T,h}(\bu,\bv) &\colonequals \int_{\Gamma_h} E_{T,h}(\bu) : E_{T,h}(\bv)\, ds_h +
                                  \int_{\Gamma_h} \bP_h \bu \cdot \bP_h \bv \, ds_h\,,\\
  b_h(\bu,q)  &\colonequals \int_{\Gamma_h} \bu \cdot \gradGh q \, ds_h\,,\quad
  \bk_h(\bu,\bv) \colonequals \eta \int_{\Gamma_h} (\bu \cdot \tilde{\bn}_h) (\bv \cdot \tilde{\bn}_h)  \, ds_h\,.
\end{align*}
The bilinear form $\bk_h(\cdot, \cdot)$ is used in a penalty approach in order to  (approximately) satisfy the condition $\bu \cdot \bn = 0$. The normal vector $\tilde{\bn}_h$, used in the penalty term $\bk_h(\cdot, \cdot)$, and the curvature tensor $\bH_h$ are approximations of the exact normal and the exact Weingarten mapping, respectively. There are several possibilities for constructing suitable approximations, e.g.,
\begin{align}\label{better_normal}
  \tilde{\bn}_h = \frac{\nabla(I_\Theta^{k+1}(\varphi))}{\Vert \nabla(I_\Theta^{k+1}(\varphi)) \Vert_2}\,, \qquad
  \bH_h = \nabla(I_\Theta^{k-1}(\bn_h))\,,
\end{align}
where $I_{\Theta}^{k} \colon L^2(\Omega_{\Theta}^{\Gamma}) \to V_{h,\Theta}^{k}$ is the parametric Oswald-type interpolation operator as defined in \cite{lehrenfeld2016high}. For these approximations we have the following error bounds:
\begin{equation} \label{estDemlow1}
  \Vert \tilde{\bn}_h - \bn \Vert_{L^{\infty}(\Omega_\Theta^\Gamma)} \lesssim h^{k+1}\,, \qquad
  \Vert \bH_h - \bH \Vert_{L^{\infty}(\Omega_\Theta^\Gamma)} \lesssim h^{k-1}\,.
\end{equation}
The reason that we introduce yet another normal approximation $\tilde{\bn}_h$  comes from error analyses for a surface vector-Laplace equation \cite{hansbo2019analysis,jankuhn2019}, which show that for obtaining optimal order estimates the normal approximation $\tilde{\bn}_h$ used in the penalty term has to be at least one order more accurate than the normal approximation $\bn_h$.

For the discretization on the surface approximation $\Gamma_h$ we  need a suitable (sufficiently accurate) extension of the data $\bbf$, which is denoted by $\bbf_h$. For $\bbf_h$ we can choose any smooth extension to the neighborhood $U_\delta$. For example, if $\bbf$ is defined on $U_\delta$ we can choose $\bbf_h = \bbf$.

\begin{remark} \label{remarkextension}
In the numerical experiments below we use the following data extension. In the setting of these experiments we prescribe an exact solution pair $(\bu,p)$ on $\Omega$ and a corresponding right hand-side $\bbf_h$ is constructed as follows. The surface  differential operators used in the Stokes problem \eqref{eq:strong}, defined on $\Gamma$, have canonical extensions to a small neighborhood of $\Gamma$. We use these extended ones and  apply the Stokes operator (defined in the neighborhood) to the prescribed $\bu$ and $p$. The resulting $\bbf$, which is defined in the neighborhood and not necessarily constant in normal direction,  is used is the numerical experiments.
\end{remark}

We now introduce a discrete version of the formulation \eqref{contform}: \\
Determine $(\bu_h, p_h) \in \bU_h \times Q_h$ such that
\begin{equation}  \label{discreteform1}
\begin{aligned}
  A_{T,h}(\bu_h,\bv_h) + b_h(\bv_h,p_h) &= (\mathbf{f}_h,\bv_h)_{L^2(\Gamma_h)} &&\text{for all } \bv_h \in \bU_h, \\
  b_h(\bu_h,q_h) - s_h(p_h,q_h) &= 0 &&\text{for all }q_h \in Q_h\,,
\end{aligned}
\end{equation}
with $A_{T,h}(\bu,\bv) \colonequals \ba_{T,h}(\bu,\bv) + \bs_h(\bu,\bv) + \bk_h(\bu,\bv)$.
Based on error analyses \cite{jankuhn2019,OlshanskiiZhiliakov2019} the stabilization parameters are chosen as
$\eta = h^{-2},\,\rho_{\bu} = h^{-1},\,\rho_p = h$. Concerning  the latter two parameters that determine the size of the velocity and pressure normal derivative stabilizations we note the following, cf.~\cite{OlshanskiiZhiliakov2019}: The stabilization for velocity is
not essential for stability of the finite element discretization method but needed (only) to control the condition number of the stiffness matrix.
The pressure stabilization  term, however, with scaling
$\rho_p  \geq c_p h$, $c_p >0$, turns out to be crucial for good (discrete inf-sup) stability properties of the finite
element discretization method.

\subsubsection{Discretization of variational formulation in stream function variable} \label{sectstream}
For the discretization of the problems \eqref{NotesAR_8}, \eqref{Variationsformulierung_u}, and \eqref{Variationsformulierung_p} we use the TraceFEM, cf.  \cite{BrandnerReusken2019}. We choose the same parametric trace finite element space for the velocity approximation as in the previous section, i.e., $\bu_h \in \bV_{h,\Theta}^{k}$. In the surface Taylor-Hood case \eqref{TaylorHood} we need  $k \geq 2$. Here, however, we allow $k \geq 1$. For the stream function approximation we also use a parametric trace finite element space, but with polynomials of one degree higher, i.e., $\psi_h \in V_{h,\Theta}^{k+1}$. We use the same stabilizations  $\bs_h(\cdot,\cdot)$, $s_h(\cdot,\cdot)$ and notations as in the previous section \ref{sec:trace_fem_th}. Note that in the geometry approximation (the parametric mapping $\Theta_h$) we use the same polynomial degree $k$ as for  the velocity approximation.

For a discrete version of \eqref{NotesAR_8} we define for $\xi,\eta \in V_{h,\Theta}^{k+1}$ the   bilinear and linear forms
\[ \begin{split}
  m_h(\xi,\eta)& \colonequals  \int_{\Gamma_h} \xi \eta \, ds_h\,, ~~
  \ell_h(\xi,\eta) \colonequals  \int_{\Gamma_h} \nablaGh \xi \cdot  \nablaGh \eta \, ds_h\,,\\
  \ell_{h,K}(\xi,\eta) & \colonequals  2 \int_{\Gamma_h} \left( 1-\tilde{K}_h \right) \nablaGh \xi \cdot  \nablaGh \eta \, ds_h\,,
~~g(\xi) \colonequals  -2 \int_{\Gamma_h} \blf_h \cdot \RotGh \xi \, ds_h.
\end{split}
\]
For the approximations of the Weingarten mapping and of the Gaussian curvature we take
\begin{equation} \label{GaussK}
  \tilde \bH_h  \colonequals  \bP_h \nabla(I_\Theta^k(\tilde \bn_h)) \bP_h, \qquad \tilde{K}_h \colonequals \frac{1}{2} \left( \tr \left( \tilde \bH_h \right)^2 - \tr \left( \tilde \bH_h^2 \right) \right),
\end{equation}
with $I_\Theta^k$ the Oswald-type interpolation also used in \eqref{better_normal}.
The approximation of the Gaussian curvature is based on the  identity $K \bP = \tr(\bH)\bH-\bH^2$ (cf.~\cite{Jankuhn1}). The following estimates hold:
\begin{equation*}
  \Vert \tilde \bH_h - \bH \Vert_{L^{\infty}(\Omega_\Theta^\Gamma)} \lesssim h^{k}, \qquad \Vert \tilde{K}_h - K \Vert_{L^{\infty}(\Omega_\Theta^\Gamma)} \lesssim h^{k}.
\end{equation*}
The discretization of \eqref{NotesAR_8}  is as follows:\\
Determine $(\phi_h,\psi_h) \in V_{h,\Theta}^{k+1}\times V_{h,\Theta}^{k+1}$ with $\int_{\Gamma_h} \psi_h \, ds_h= 0$,  such that
\begin{equation}\label{NotesAR_10}
\begin{aligned}
  m_h(\phi_h,\eta_h) + \ell_h(\psi_h, \eta_h) + s_h(\psi_h, \eta_h) &= 0 &&\text{for all } \eta_h \in V_{h,\Theta}^{k+1}\,,\\
  \ell_h(\phi_h,\xi_h) - \ell_{h,K}(\psi_h, \xi_h) + s_h(\phi_h, \xi_h) &= g_h(\xi_h) &&\text{for all }\xi_h \in V_{h,\Theta}^{k+1}\,.
\end{aligned}
\end{equation}
Based on the analysis in \cite{BrandnerReusken2019}, the parameters in both stabilizations $s_h(\cdot,\cdot)$ in \eqref{NotesAR_10} are set to $\rho_p = h$. For the discretization of the velocity reconstruction we introduce the bilinear form
$
  \bm_h(\bu, \bv) \colonequals \int_{\Gamma_h} \bu \cdot \bv \, ds_h $
and define the discrete problem: Determine $\bu_h \in \bV_{h,\Theta}^{k}$ such that
\begin{equation}\label{Variationsformulierung_u_diskret}
  \bm_h(\bu_h, \bv_h) + \bs_h(\bu_h,\bv_h) = \int_{\Gamma_h} (\tilde{\bn}_h \times \nablaGh \psi_h)\cdot \bv_h \,ds_h \quad \forall~\bv_h \in \bV_{h,\Theta}^{k}\,,
\end{equation}
with the given discrete solution $\psi_h$ of \eqref{NotesAR_10} and the normal approximation $\tilde{\bn}_h$ as in \eqref{better_normal}. In the stabilization bilinear form $\bs_h(\cdot,\cdot)$ we use the parameter $\rho_{\bu} = h$, based on the analysis in \cite{BrandnerReusken2019}.
Concerning the reconstruction of the pressure we consider the following discrete variational formulation of \eqref{Variationsformulierung_p}: Determine $p_h \in V_{h,\Theta}^{k}$ with $\int_{\Gamma_h} p_h \, ds_h=0$,  such that
\begin{equation}\label{Variationsformulierung_p_diskret}
  \ell_h(p_h , \xi_h) + s_h(p_h,\xi_h) = \int_{\Gamma_h} \left( \tilde{K}_h \RotGh \psi_h + \blf_h \right) \cdot \nablaGh \xi_h \,ds_h \quad \forall~\xi_h \in V_{h,\Theta}^{k}\,.
\end{equation}
We use the parameter choice $\rho_p = h$ in the stabilization $s_h(\cdot,\cdot)$, cf. Table~\ref{tabpar}.

\subsection{Surface finite element method}
The surface finite element discretization combines the general piecewise flat surface discretization for vector-valued surface partial differential equations \cite{Nestleretal_arxiv_2018} with the higher order surface approximations considered for scalar-valued surface partial differential equations \cite{demlow2009higher}, which requires some additional handling of the tangential constraint. The stream function formulation, on the other hand, is a straightforward extension of the piecewise flat surface discretization \cite{nitschke2012finite} to curved geometries.

\subsubsection{Discretization of variational formulation in $\bu$-$p$ variables} \label{sectupS}
Similar to the TraceFEM discretization, we introduce the $\vect P_k$-$P_{k-1}$ pair of surface Taylor-Hood elements based on the function spaces $\tilde \bV_{h,\pi_h}^{k}$ and $\tilde V_{h,\pi_h}^{k}$:
\begin{equation} \label{TaylorHoodS}
  \tilde \bU_h \colonequals \tilde \bV_{h,\pi_h}^{k},\quad
  \tilde Q_h \colonequals \tilde V_{h,\pi_h}^{k-1}\cap L_0^2(\Gamma_h),\quad k\geq 2\,.
\end{equation}
Using the same discrete bilinear forms $\ba_{T,h}(\cdot,\cdot)$, $b_h(\cdot,\cdot)$, and $\bk_h(\cdot,\cdot)$ as for the TraceFEM discretization, but applied to functions from $\tilde \bU_h$ and $\tilde Q_h$, we can directly formulate the discrete problem. For the approximation $\tilde{\bn}_h$ of $\bn$ in the penalty term $\bk_h(\cdot,\cdot)$, we use a local Lagrange interpolation of the exact surface normal $\bn(\xx) = \nabla\varphi(\xx)/\|\nabla\varphi(\xx)\|$ for $\xx\in\Gamma$, i.e.,
\begin{equation}\label{eq:better_normal_approximation}
  \tilde{\bn}_h = \frac{\tilde{\bN}_h}{\|\tilde{\bN}_h\|},\;\text{ with }\;\tilde{\bN}_h\circ \bar{\pi}_{\hat{S},h} = I^k(\bn\circ\pi\circ \bar{\pi}_{\hat{S},h})\,,
\end{equation}
resulting in an approximation with $\|\tilde{\bn}_h - \bn\circ \pi\|_{L^\infty(\Gamma_h)} \lesssim h^{k+1}$ that follows from standard interpolation error estimates.

We now introduce a discrete version of the formulation \eqref{contform}:
Determine $(\bu_h, p_h) \in \tilde{\bU}_h \times \tilde{Q}_h$ such that
\begin{equation}  \label{discreteform2}
\begin{aligned}
  A_{T,h}(\bu_h,\bv_h) + b_h(\bv_h,p_h) &= (\mathbf{f}_h,\bv_h)_{L^2(\Gamma_h)} &&\text{for all } \bv_h \in \tilde{\bU}_h \\
  b_h(\bu_h,q_h) &= 0 &&\text{for all }q_h \in \tilde{Q}_h,
\end{aligned}
\end{equation}
with $A_{T,h}(\bu,\bv) \colonequals \ba_{T,h}(\bu,\bv) + \bk_h(\bu,\bv)$.
Based on error analysis in \cite{hansbo2019analysis} the penalization parameter is chosen as $\eta = h^{-2}$. We use the same data extension $\mathbf{f}_h=\mathbf{f}$ on $\Gamma_h$ as explained in Remark \ref{remarkextension}. Note that compared to the TraceFEM discretization no additional stabilization terms, such as  $\bs_h(\cdot,\cdot)$, are needed.

\subsubsection{Discretization of variational formulation in stream function variable} \label{sectstreamS}
The SFEM discretization of \eqref{NotesAR_8}, \eqref{Variationsformulierung_u}, and \eqref{Variationsformulierung_p} is analogous to the TraceFEM discretization but without the stabilization terms.

We use the scalar finite element approximation $\psi_h\in\tilde{V}_{h,\pi_h}^{k+1}$ for the stream function and $\phi_h\in\tilde{V}_{h,\pi_h}^{k+1}$ for the vorticity function of one order higher than the corresponding velocity function $\bu_h\in\tilde{\bV}_{h,\pi_h}^{k}$, but on the surface approximation $\Gamma_h$ that is constructed using polynomials of degree $k$. The bilinear forms $m_h(\cdot,\cdot)$ and $\ell_h(\cdot,\cdot)$ are defined  as in the TraceFEM discretization. The bilinear form $\ell_{h,H}(\cdot,\cdot)$ for the SFEM discretization requires a higher order accurate curvature tensor  $\tilde{\bH}_h=\nabla_S\tilde{\bn}_h$, given by
\begin{equation*}
  \tilde{\bH}_h = \nabla_{S}\tilde{\bN}_h\, \frac{(I - \tilde{\bn}_h\tilde{\bn}_h^T)}{\|\tilde{\bN}_h\|}\,,
\end{equation*}
locally, inside each element $S\in\mathcal{S}_h$, with $\tilde{\bn}_h$ and $\tilde{\bN}_h$ from \eqref{eq:better_normal_approximation}.
This leads to the following discrete version of the equation \eqref{NotesAR_8}:
Determine $(\phi_h, \psi_h) \in\tilde{V}_{h,\pi_h}^{k+1} \times \tilde{V}_{h,\pi_h}^{k+1}$ with $\int_{\Gamma_h} \psi_h \, ds_h= 0$,  such that
\begin{equation}\label{discreteform_sf2}
\begin{aligned}
  m_h(\phi_h,\eta_h) + \ell_h(\psi_h, \eta_h) &= 0 &&\text{for all } \eta_h \in\tilde{V}_{h,\pi_h}^{k+1}\,,\\
  \ell_h(\phi_h,\xi_h) - \ell_{h,K}(\psi_h, \xi_h)  &= g_h(\xi_h) &&\text{for all } \xi_h \in\tilde{V}_{h,\pi_h}^{k+1}\,.
\end{aligned}
\end{equation}

With $\bm_h(\cdot,\cdot)$ defined as for the TraceFEM discretization, we obtain for the discretization of the velocity reconstruction the discrete problem:
Determine $\bu_h \in\tilde{\bV}_{h,\pi_h}^{k}$ such that
\begin{equation}\label{discreteform_u2}
  \bm_h(\bu_h , \bv_h) = \int_{\Gamma_h} (\tilde{\bn}_h \times \nablaGh \psi_h)\cdot \bv_h \,ds_h \quad \forall~\bv_h \in\tilde{\bV}_{h,\pi_h}^{k}\,,
\end{equation}
with the given discrete solution $\psi_h$ of \eqref{discreteform_sf2} and the normal approximation $\tilde{\bn}_h$ as before.
Concerning the reconstruction of the pressure we consider the following discrete variational formulation of \eqref{Variationsformulierung_p}:
Determine $p_h \in\tilde{V}_{h,\pi_h}^{k}$ with $\int_{\Gamma_h} p_h \, ds_h=0$,  such that
\begin{equation}\label{discreteform_p2}
  \ell_h(p_h , \xi_h) = \int_{\Gamma_h} \left( \tilde{K}_h \RotGh \psi_h + \blf_h \right) \cdot \nablaGh \xi_h \,ds_h \quad \forall~\xi_h \in\tilde{V}_{h,\pi_h}^{k}.
\end{equation}
using $\tilde{K}_h$ given by $ \frac{1}{2}\big(\tr(\tilde\bH_h)^2 - \tr(\tilde\bH_h^2) \big)$.
Again, no additional stabilization terms are needed.

\subsection{Comparison of formulations and methods}

\subsubsection{Comparison of velocity-pressure and stream function formulations} \label{sectdisc}
Relations between the velocity-pressure formulation and the stream function formulation for the Stokes equations in  Euclidean space are treated in e.g. \cite{GR}. In Sections~\ref{sec:trace_fem_th} and \ref{sectstream} additional properties associated with the \emph{surface} Stokes equations are discussed. While in the Taylor-Hood formulation a penalty approach is used to (weakly) enforce the tangential condition $\bu \cdot \bn=0$ we do not need an additional Lagrange multiplier or penalty approach to enforce the tangential condition in the stream function formulation. If the velocity is reconstructed based on  the relation $\bu^\ast_T= \bn \times \gradG \psi^\ast$, tangency is automatically fulfilled. In both formulations one needs curvature information. In the Taylor-Hood method (both for TraceFEM and SFEM) an approximation $\bH_h$ of the exact Weingarten mapping is needed in the bilinear form $a_{T,h}(\cdot,\cdot)$. This approximation should have accuracy (at least) $h^{k-1}$, cf. \eqref{estDemlow}, \eqref{estDemlow1}. In the stream function approach we need an approximation of the Gaussian curvature $K$, which is determined based on an approximation of the Weingarten mapping as in \eqref{GaussK}. This Gaussian curvature approximation should have accuracy (at least) $h^k$.  In both formulations a one order more accurate normal approximation (denoted by $\tilde \bn_h$ above) is used: in the Taylor-Hood formulation in the penalty bilinear form $\bk_h(\cdot,\cdot)$ and in the stream function formulation in the reconstruction of the velocity, cf. \eqref{Variationsformulierung_u_diskret}, \eqref{discreteform_u2}.

The assumption that $\Gamma$ is simply connected is crucial for the stream function formulation. Note that the Helmholtz-Hodge decomposition, which provides the mathematical basis for the stream function formulation, not only splits $\bu_T$ into curl-free and divergence-free components, but might also contain non-trivial harmonic vector fields -- vector fields which are curl- and divergence-free. As these vector fields cannot be described by the stream function formulation, the approach is only applicable for surfaces, where harmonic vector fields are trivial, which are only simply-connected surfaces, see also \cite{Nitschke_DEC_2017,pamm2021} for numerical comparisons. Such a restriction does not exist for the Taylor-Hood formulation.

\subsubsection{Comparison of TraceFEM and SFEM} \label{sectC}
We compare  complexity of both discretization methods for  the two formulations (mixed and stream function) in terms of the number of degrees of freedom (DOFs) for representing the discrete solution. Clearly, the number of unknowns  depends on the underlying mesh. As a theoretical model case we consider a \emph{flat} surface  with a structured (uniform) triangulation and derive formulas for  the number of unknowns for that case. We then  perform a numerical experiment to test how well these formulas predict the number of unknowns for the Stokes problem on a sphere discretized using quasi-uniform outer- (TraceFEM) or surface-triangulations (SFEM).

For TraceFEM we consider the following structured case. Let $\Omega=[-1,1]^3$ and $\Gamma = \{ \xx=(x,y,z) \in \mathbb{R}^3 \mid z = \frac{1}{3} \}$. We assume periodic boundary conditions on the faces of the cube. The initial triangulation consists of $4^3$ equal sized sub-cubes, where each of these is subdivided into $6$ equal sized tetrahedra. For refinement, each sub-cube is repeatedly divided into $8$ sub-cubes.
For SFEM the procedure for the construction is similar, but we just consider the 2d level-set domain $\Gamma$. The initial triangulation consists of $4^2$ equal sized quads, where each of these is subdivided into 2 triangles. For refinement, each sub-quad is repeatedly divided into 4 sub-quads.

On these structured meshes we compare the number of unknowns for both formulations \eqref{discreteform1} and \eqref{NotesAR_10} as a function of the number of tetrahedra cut by the surface $\Gamma$ (TraceFEM) and triangular surface-grid elements (SFEM), which we denote by $n$, and the degree $k$ of the finite elements. We start with the Taylor-Hood formulations \eqref{discreteform1} and \eqref{discreteform2}. For the number of unknowns for the $\vect P_k$-$P_{k-1}$ pair of parametric Taylor-Hood elements one can derive the formulas
\begin{equation} \label{nr1}\begin{split}
  TH^\text{TraceFEM}(n,k) &= 3 \Big\lceil \frac{n}{6} \Big\rceil (k+1)k^2 + \Big\lceil \frac{n}{6} \Big\rceil k(k-1)^2\,, \\
  TH^\text{SFEM}(n,k) &= 3 \Big\lceil \frac{n}{2} \Big\rceil k^2 + \Big\lceil \frac{n}{2} \Big\rceil (k-1)^2\,.
\end{split}
\end{equation}

The first summand is the number of unknowns corresponding to $\vect P_k$ (velocity) and the second summand is the number of unknowns corresponding to $P_{k-1}$ (pressure). For the stream function formulations \eqref{NotesAR_10} and \eqref{discreteform_sf2} the number of unknowns for the coupled stream function/vorticity system discretized with $P_{k+1}$ finite elements both for the stream function and the vorticity is given by
\begin{equation} \label{nr2}\begin{split}
  SF^\text{TraceFEM}(n,k) &= 2 \Big\lceil \frac{n}{6} \Big\rceil (k+2)(k+1)^2\,, \\
  SF^\text{SFEM}(n,k) &= 2 \Big\lceil \frac{n}{2} \Big\rceil (k+1)^2\,.
\end{split}
\end{equation}
Adding the number of unknowns for reconstructing the velocity in \eqref{Variationsformulierung_u_diskret} and \eqref{discreteform_u2} and the pressure in \eqref{Variationsformulierung_p_diskret} and \eqref{discreteform_p2}, using finite elements of degree $k$, we get
\begin{equation} \label{nr3}\begin{split}
  SF^\text{TraceFEM}_{total}(n,k) &= 2 \Big\lceil \frac{n}{6} \Big\rceil (k+2)(k+1)^2 + (3 + 1) \Big\lceil \frac{n}{6} \Big\rceil (k+1)k^2\,, \\
  SF^\text{SFEM}_{total}(n,k) &= 2 \Big\lceil \frac{n}{2} \Big\rceil (k+1)^2 + (3 + 1) \Big\lceil \frac{n}{2} \Big\rceil k^2\,.
\end{split}
\end{equation}
Differences between TraceFEM and SFEM are the constant factor $\lceil \frac{n}{6}\rceil$ vs. $\lceil\frac{n}{2}\rceil$, which relates to the number of simplices in the corresponding cube/quad elements, and the missing third dimension for the SFEM discretization. The latter yields one polynomial degree lower dependency on $k$ in SFEM discretizations compared to the TraceFEM discretizations.

\begin{figure}[ht]
  \centering
\def\nnn{102242}
\def\nn{131072}
\begin{subfigure}[b]{0.49\textwidth}
  \scalebox{0.80}{\begin{tikzpicture}
    \begin{axis}[xlabel={$k$}, width=1.3\linewidth, height=\linewidth,
                 legend entries={$TH$,$SF$,$SF_{total}$},
                 legend pos=north west]
      \addplot[blue, domain=2:6,samples=201,]{3*ceil(\nnn/6)*(x+1)*x^2 + ceil(\nnn/6)*x*(x-1)^2};
      \addplot[red, domain=2:6,samples=201,]{2*ceil(\nnn/6)*(x+2)*(x+1)^2};
      \addplot[brown!60!black, domain=2:6,samples=201,]{2*ceil(\nnn/6)*(x+2)*(x+1)^2 + 3*ceil(\nn/6)*(x+1)*x^2 + ceil(\nnn/6)*(x+1)*x^2};
      \addplot[blue,dashed, domain=2:6,samples=201,]{3*ceil(\nn/2)*x^2 + ceil(\nn/2)*(x-1)^2};
      \addplot[red,dashed, domain=2:6,samples=201,]{2*ceil(\nn/2)*(x+1)^2};
      \addplot[brown!60!black,dashed, domain=2:6,samples=201,]{2*ceil(\nn/2)*(x+1)^2 + 3*ceil(\nn/2)*x^2 + ceil(\nn/2)*x^2};
    \end{axis}
  \end{tikzpicture}}
  \caption{Structured plane}
\end{subfigure}%
\begin{subfigure}[b]{0.49\textwidth}
  \scalebox{0.80}{\begin{tikzpicture}
    \begin{axis}[xlabel={$k$}, width=1.3\linewidth, height=\linewidth,
                 legend entries={$TH$,$SF$,$SF_{total}$},
                 legend pos=north west]
      \addplot[blue,mark=square*] table [x=k, y=dofTH] {comp.dat};
      \addplot[red,mark=triangle*] table [x=k, y=dofSF] {comp.dat};
      \addplot[brown!60!black,mark=*] table [x=k, y=dofSFtotal] {comp.dat};
      \addplot[blue,dashed,mark=square*,mark options={solid}] table [x=k, y=dofTH] {SFEM_comp.dat};
      \addplot[red,dashed,mark=triangle*,mark options={solid}] table [x=k, y=dofSF] {SFEM_comp.dat};
      \addplot[brown!60!black,dashed,mark=*,mark options={solid}] table [x=k, y=dofSFtotal] {SFEM_comp.dat};
    \end{axis}
  \end{tikzpicture}}
  \caption{Unstructured sphere}
\end{subfigure}
   \caption{Plot of formulas \eqref{nr1}, \eqref{nr2}, and \eqref{nr3}, for number of unknowns on a structured grid (left) and number of unknowns for  the unit sphere and an unstructured triangulation consisting of $n= 100474$ tetrahedra cut by the surface for TraceFEM and $n=131072$ triangles for SFEM, respectively (right). Solid lines correspond to TraceFEM discretizations whereas dashed lines correspond to SFEM discretizations (Color coding available online).}\label{fig:compunknowns}
\end{figure}
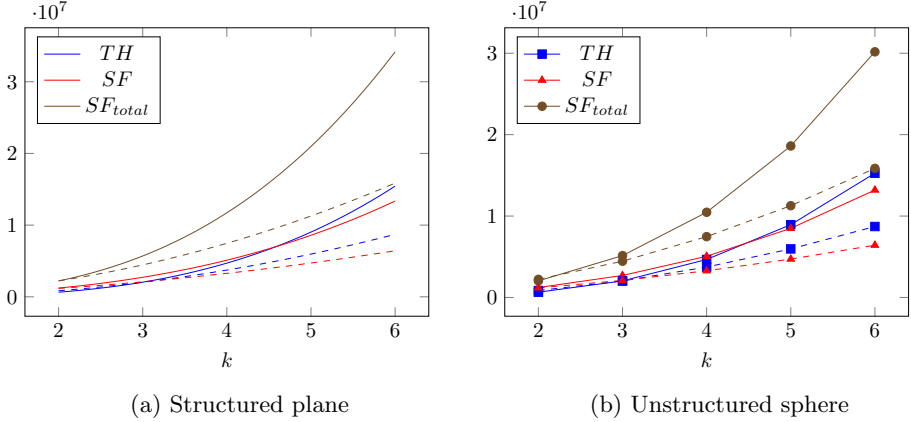

In Figure \ref{fig:compunknowns} we illustrate these formulas. On the left-hand side of Figure \ref{fig:compunknowns} we plotted the formulas \eqref{nr1}, \eqref{nr2}, and \eqref{nr3} with a fixed $n$ extracted from the surface grid experiments that are plotted on the right-hand side. We use solid lines for the TraceFEM discretizations and dashed lines for the SFEM discretizations.
On the right-hand side of that figure we plotted as a comparison the number of unknowns for the case that $\Gamma$ is the unit sphere and the computational grids are unstructured, consisting of $n=100474 $ tetrahedra (TraceFEM) and $n=131072$ triangles (SFEM). We observe that the corresponding curves in the two figures are very close.

From these formulas it follows that $TH(n,k) < SF(n,k)$ for $k\leq 4$ (TraceFEM) or $k\leq 3$ (SFEM) and $TH(n,k) > SF(n,k)$ for  $k\geq 5$ (TraceFEM) or $k\geq 4$ (SFEM). Thus, for low polynomial degree in the finite elements, the Taylor-Hood formulation has an advantage in terms of number of degrees of freedom over the stream function formulation. The difference, however,  between the number of unknowns in   these two formulations is relatively small.
If we include the number of unknowns for reconstructing the velocity and pressure we obtain $TH(n,k) < SF_{total}(n,k)$ for $k\geq 2$. This is due to the fact that the reconstruction of the velocity in the second formulation is performed in the same finite element space as the one used for velocity in the first formulation. This choice is necessary to obtain the same convergence order of convergence, see below.

Note that we only measure the number of unknowns and \emph{not} the computational costs of an iterative (or sparse direct) method for solving the resulting linear systems.

\paragraph{System complexity}
As already indicated by the formal analysis of number of unknowns in the discretizations on a plane or sphere, the TraceFEM and SFEM have different complexities. This is due to the fact that the TraceFEM is based on finite element spaces on a strip of 3D tetrahedral elements whereas the SFEM is based on spaces on 2D triangular elements. This difference does not only imply a higher number of unknowns for the TraceFEM compared to SFEM, but also increases the number of non-zeros (nnz) in the resulting linear systems. For a computational grid used in Section \ref {sec:numerical_experiments}, we have listed complexity data for TraceFEM and SFEM, for the case $k=3$, in the Table \ref{tab:comparison_data}.

We note that  the two methods have a different mesh size parameter $h$. In TraceFEM it is natural to use a mesh size parameter $h$ corresponding to the  bulk tetrahedral mesh. In SFEM the mesh size $h$ is as usual  the longest edge length in the surface triangulation.
For reasons of comparison, for TraceFEM we also determined the surface grid size, denoted by $h_\Gamma$, defined as the average over the longest edges in the triangular faces of the surface cut through the tetrahedral computational grid elements.

\begin{table}[ht]
  \begin{center}
  \begin{tabular}{l|l|l|l|l|l|l}
    \toprule
    ref.  & grid         & grid     &  \multicolumn{2}{c|}{TH}      & \multicolumn{2}{c}{SF} \\
    level	& size         & elements & nnz & DOFs                    & nnz & DOFs \\
    \hline
    TraceFEM &  $h_\Gamma$ \\
    \hline
    $0$   & $0.288439$   & $334$    & $5.03 \cdot 10^5$ & $6843$    & $5.99 \cdot 10^5$ & $9054$ \\
    $2$   & $0.078671$   & $6250$   & $9.62 \cdot 10^6$ & $126724$  & $1.14 \cdot 10^7$ & $168300$ \\
    $4$   & $0.019486$   & $102640$  & $1.58 \cdot 10^8$ & $2077692$ & $1.87 \cdot 10^8$ & $2760602$ \\
    \hline
    SFEM &  $h$ \\
    \hline
    $0$   & $0.282815$ & $516$    & $4.88\cdot 10^5$ & $8006$    & $5.46\cdot 10^5$ & $8260$  \\
    $2$   & $0.075342$ & $8256$   & $7.81\cdot 10^6$ & $127976$  & $8.73\cdot 10^6$ & $132100$ \\
    $4$   & $0.019885$ & $132096$ & $1.25\cdot 10^8$ & $2047496$ & $1.40\cdot 10^8$ & $2113540$ \\
    \bottomrule
  \end{tabular}
  \end{center}
  \caption{Complexity of TraceFEM (top) and SFEM (bottom) grids: (Surface) grid size, number of elements in the computational grid,  number of non-zeros (nnz) and number of degrees of freedom (DOFs) in the two formulations for $k=3$ on the ``biconcave shape'' grid.}\label{tab:comparison_data}
\end{table}

\subsubsection{Comparison of numerical parameters}
We summarize the  parameter settings used in the TraceFEM of Section~\ref{sec:trace_fem_th} (TraceFEM TH), the TraceFEM of Section~\ref{sectstream} (TraceFEM SF) and the surface FEM of Section~\ref{sectupS} (SFEM TH) in Table~\ref{tabpar}.

\begin{table}[ht]
  \centering
  {\renewcommand{\arraystretch}{1.1}
  \begin{tabular}{r|c|c|c|c}
    \toprule
      & $\rho_{\mathbf{u}}$ & $\rho_{p}$ & $\rho_{\psi}$ & $\eta$ \\
    \hline
    TraceFEM TH & $h^{-1}$ & $h$ &   & $h^{-2}$ \\
    SFEM TH &  &  &  & $h^{-2}$ \\
    TraceFEM SF & $h$ & $h$ & $h$ & \\
    \bottomrule
  \end{tabular}
  \caption{Choices for the normal derivative volume stabilization parameters $\rho_{\mathbf{u}}$, $\rho_{p}$, $\rho_{\psi}$ and the penalty parameter $\eta$ in the different methods.}
  \label{tabpar}}
\end{table}

\section{Numerical experiments}\label{sec:numerical_experiments}
The implementation of the TraceFEM discretizations is done with Netgen/NGSolve and ngsxfem \cite{Schoeberl2014,ngsxfem}. We use an unstructured tetrahedral triangulation of $\Omega \colonequals [-\frac{5}{3},\frac{5}{3}]^3$ with starting mesh size $h=0.5$ in the three-dimensional grid. The mesh is locally refined using a marked-edge bisection method for the surface intersected tetrahedra \cite{Schoeberl1997}. The average over the longest edges in the triangular faces of the surface cut through the initial tetrahedral computational grid elements is about $0.29$. The code is available in \cite{Brandner2021CodeTraceFEMSF,Jankuhn2021CodeTraceFEMTH}.

The implementation of the SFEM discretizations is done with AMDiS/Dune \cite{AMDiS,Witkowskietal_ACM_2015,Dune,DuneBook}. We use an unstructured curved triangular grid \cite{praetorius2020dune,SaKoScFl2017} with initial mesh size $h\approx 0.28$ using a quartering (red) refinement of all triangles. The code is available in \cite{Praetorius2021CodeSFEM}.

\subsection{Experiment on a biconcave shape}
\begin{figure}[ht!]
  \centering
  \begin{subfigure}[b]{0.31\textwidth}
  \hspace*{-0.5cm}
      \centering
      \includegraphics[width=1.3\textwidth]{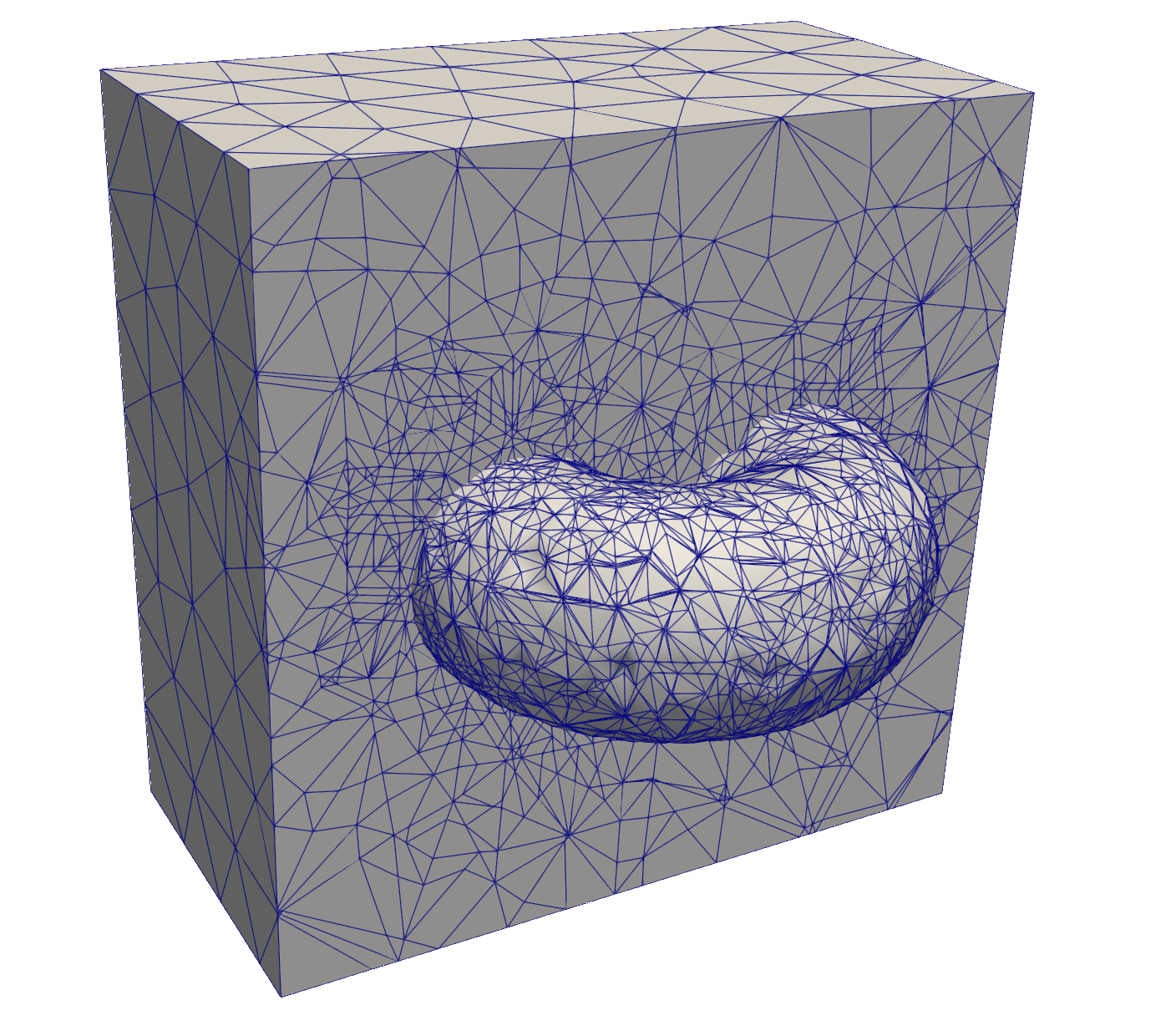}
      \caption{TraceFEM triangulation}
  \end{subfigure}
  \hfill
  \begin{subfigure}[b]{0.32\textwidth}
      \centering
      \includegraphics[width=\textwidth]{figure_2_1.png}
      \caption{SFEM triangulation}
  \end{subfigure}
  \hfill
  \begin{subfigure}[b]{0.32\textwidth}
      \centering
      \includegraphics[width=\textwidth]{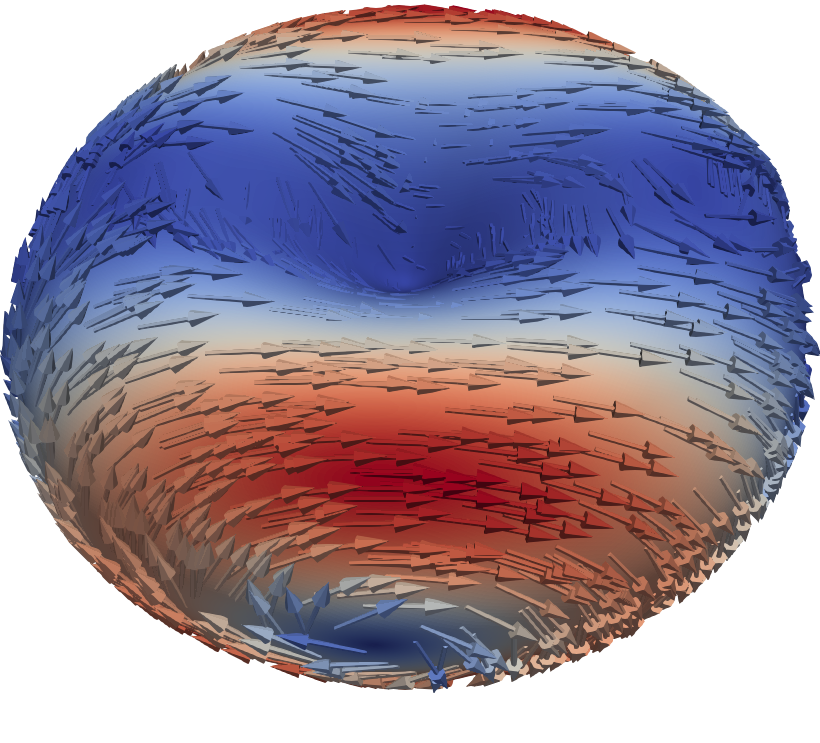}
      \caption{Velocity solution $\bu$}
  \end{subfigure}
  \caption{Geometry ``biconcave shape'' with triangulations and solution \eqref{eq:exact_solutions}. In the TraceFEM triangulation a cut through the bulk mesh and cut surface elements of the levelset are shown. Glyphs in the solution represent the velocity, the coloring from blue to red the velocity magnitude (Color coding available online).}\label{fig:biconcave-shape}
\end{figure}

We consider a surface $\Gamma$ called ``biconcave shape'', cf. Fig.~\ref{fig:biconcave-shape}, which has high curvature variations and is defined implicitly as the zero-level set of a function $\varphi(\xx)$:
\begin{equation}\label{eq:biconcave-shape_levelset}
  \Gamma \colonequals \{\xx=(x,y,z) \in \Bbb{R}^3\mid \varphi(\xx) \colonequals (d^2 + x^2 + y^2 + z^2)^3 - 8 d^2 (y^2 + z^2) - c^4 = 0 \}
\end{equation}
with $c=0.95$ and $d=0.96$. The smooth solutions of the surface Stokes problem \eqref{eq:strong} are prescribed by
\begin{equation}\label{eq:exact_solutions}
  p   \colonequals x^3 + xyz\,, \quad
  \bu \colonequals \RotG \psi\,,\text{ and }\,
  \psi\colonequals x^2y - 5z^3.
\end{equation}

Note that $\bu$ is tangential and divergence-free. The solutions are not extended constantly along normals but by their definitions \eqref{eq:exact_solutions}. Using MAPLE we calculate the corresponding right-hand side $\blf$ as described in Remark \ref{remarkextension}. The problem setting is illustrated in Fig.~\ref{fig:biconcave-shape}.

\begin{figure}[ht]
  \centering
\begin{subfigure}[b]{0.49\textwidth}
  \scalebox{0.80}{\begin{tikzpicture}
  \def\vara{5}
  \def\varb{120}
  \begin{loglogaxis}[ xlabel={Mesh size $h$},domain={0.015625:0.29},
    width=1.3\linewidth, height=\linewidth,
    legend style={at={(0.99,0.015)},anchor=south east,legend columns=2,draw=none, fill=none},
    legend cell align=left,
    cycle list name=mark list ]
    \addplot[mark=triangle] table[x=h, y=L2(U-Uh)____] {TH_k2_apple.dat};
    \addplot[mark=square] table[x=h, y=L2(U-Uh)____] {SF_k2_apple.dat};
    \addplot[mark=triangle*] table[x=h, y=L2(U-Uh)____] {TH_k3_apple.dat};
    \addplot[mark=square*] table[x=h, y=L2(U-Uh)____] {SF_k3_apple.dat};
    \addplot[dotted,line width=0.75pt] {1000*x^3};
    \addplot[dashed,line width=0.75pt] {100*x^4};
    \legend{
      \small TH,
      \small SF ($k=2$),
      \small TH,
      \small SF ($k=3$),
      \small $\mathcal{O}(h^3)$,
      \small $\mathcal{O}(h^4)$
    }
  \end{loglogaxis}
  \node[anchor=north west] at (0.1,4.6) {(a) TraceFEM};
  \end{tikzpicture}}
\end{subfigure}
\begin{subfigure}[b]{0.49\textwidth}
  \scalebox{0.80}{\begin{tikzpicture}
    \begin{loglogaxis}[ xlabel={Mesh size $h$},domain={0.015625:0.29},
      width=1.3\linewidth, height=\linewidth,
      legend style={at={(0.99,0.015)},anchor=south east,legend columns=2,draw=none, fill=none},
      legend cell align=left, cycle list name=mark list ]
      \addplot[mark=triangle] table[x=h, y=L2(U-Uh)____] {SFEM_TH_k2_apple.dat};
      \addplot[mark=square] table[x=h, y=L2(U-Uh)____] {SFEM_SF_k2_apple.dat};
      \addplot[mark=triangle*] table[x=h, y=L2(U-Uh)____] {SFEM_TH_k3_apple.dat};
      \addplot[mark=square*] table[x=h, y=L2(U-Uh)____] {SFEM_SF_k3_apple.dat};
      \addplot[dotted,line width=0.75pt] {50*x^3};
      \addplot[dashed,line width=0.75pt] {0.5*x^4};
      \legend{
        \small TH,
        \small SF ($k=2$),
        \small TH,
        \small SF ($k=3$),
        \small $\mathcal{O}(h^3)$,
        \small $\mathcal{O}(h^4)$
      }
    \end{loglogaxis}
    \node[anchor=north west] at (0.1,4.6) {(b) SFEM};
  \end{tikzpicture}}
\end{subfigure}%
\caption{$\|\bu_h-\bu\|_{L^2(\Gamma_h)}$ errors for $k=2$ and $k=3$ for the different methods.}\label{fig:biconcave-shape_U_L2}
\end{figure}

\begin{figure}[ht]
  \centering
\begin{subfigure}[b]{0.49\textwidth}
\scalebox{0.80}{\begin{tikzpicture}
  \begin{loglogaxis}[xlabel={Mesh size $h$},domain={0.015625:0.29},
    width=1.3\linewidth, height=\linewidth,
    legend style={at={(0.99,0.015)},anchor=south east,legend columns=2,draw=none, fill=none},
    legend cell align=left, cycle list name=mark list ]
    \addplot[mark=triangle] table[x=h, y=H1_sem(U-Uh)] {TH_k2_apple.dat};
    \addplot[mark=square] table[x=h, y=H1_sem(U-Uh)] {SF_k2_apple.dat};
    \addplot[mark=triangle*] table[x=h, y=H1_sem(U-Uh)] {TH_k3_apple.dat};
    \addplot[mark=square*] table[x=h, y=H1_sem(U-Uh)] {SF_k3_apple.dat};
    \addplot[dotted,line width=0.75pt] {800*x^2};
    \addplot[dashed,line width=0.75pt] {300*x^3};
    \legend{
      \small TH,
      \small SF ($k=2$),
      \small TH,
      \small SF ($k=3$),
      \small $\mathcal{O}(h^2)$,
      \small $\mathcal{O}(h^3)$
    }
  \end{loglogaxis}
  \node[anchor=north west] at (0.1,4.6) {(a) TraceFEM};
  \end{tikzpicture}}
\end{subfigure}
\begin{subfigure}[b]{0.49\textwidth}
  \scalebox{0.80}{\begin{tikzpicture}
  \begin{loglogaxis}[xlabel={Mesh size $h$},domain={0.015625:0.29},
    width=1.3\linewidth, height=\linewidth,
    legend style={at={(0.99,0.015)},anchor=south east,legend columns=2,draw=none, fill=none},
    legend cell align=left, cycle list name=mark list ]
    \addplot[mark=triangle] table[x=h, y=H1_sem(U-Uh)] {SFEM_TH_k2_apple.dat};
    \addplot[mark=square] table[x=h, y=H1_sem(U-Uh)] {SFEM_SF_k2_apple.dat};
    \addplot[mark=triangle*] table[x=h, y=H1_sem(U-Uh)] {SFEM_TH_k3_apple.dat};
    \addplot[mark=square*] table[x=h, y=H1_sem(U-Uh)] {SFEM_SF_k3_apple.dat};
    \addplot[dotted,line width=0.75pt] {30*x^2};
    \addplot[dashed,line width=0.75pt] {6*x^3};
    \legend{
      \small TH,
      \small SF ($k=2$),
      \small TH,
      \small SF ($k=3$),
      \small $\mathcal{O}(h^2)$,
      \small $\mathcal{O}(h^3)$
    }
  \end{loglogaxis}
  \node[anchor=north west] at (0.1,4.6) {(b) SFEM};
  \end{tikzpicture}}
\end{subfigure}%
\caption{$\|\nablaGh \left( \bu_h-\bu \right) \|_{L^2(\Gamma_h)}$ errors for $k=2$ and $k=3$ for the different methods.}\label{fig:biconcave-shape_U_H1}
\end{figure}

In Figures \ref{fig:biconcave-shape_U_L2} and \ref{fig:biconcave-shape_U_H1} we have plotted the discretization errors of the discrete solution $\bu_h$ vs. the exact solution $\bu$ in the $L^2$-norm and $H^1$-norm. We clearly see that the error in the $L^2$-norm is one order higher than the error in the $H^1$-norm, as expected.  \emph{For both methods (TraceFEM and SFEM) and for both formulations (mixed Taylor-Hood and stream function) we observe optimal orders $k+1$ and $k$ of the $L^2$-error and $H^1$-error convergence}, respectively, consistent with analyses presented in \cite{jankuhn2019,BrandnerReusken2019,hansbo2019analysis}.

\begin{figure}[ht]
  \centering
\begin{subfigure}[b]{0.49\textwidth}
  \scalebox{0.80}{\begin{tikzpicture}
    \begin{loglogaxis}[ xlabel={Mesh size $h$},domain={0.015625:0.29},
      width=1.3\linewidth, height=\linewidth,
      legend style={at={(0.99,0.015)},anchor=south east,legend columns=3,draw=none, fill=none},
      legend cell align=left, cycle list name=mark list ]
      \addlegendimage{empty legend}
      \addplot[mark=triangle] table[x=h, y=L2(P-Ph)____] {TH_k2_apple.dat};
      \addplot[mark=square] table[x=h, y=L2(P-Ph)____] {SF_k2_apple.dat};

      \addlegendimage{empty legend}
      \addplot[mark=triangle*] table[x=h, y=L2(P-Ph)____] {TH_k3_apple.dat};
      \addplot[mark=square*] table[x=h, y=L2(P-Ph)____] {SF_k3_apple.dat};

      \addplot[loosely dotted,line width=0.75pt] {50*x^2};
      \addplot[dotted,line width=0.75pt] {100*x^3};
      \addplot[dashed,line width=0.75pt] {80*x^4};
      \legend{
        ~,\small TH, \small SF ($k=2$),
        ~,\small TH, \small SF ($k=3$),
        \small $\mathcal{O}(h^2)$,
        \small $\mathcal{O}(h^3)$,
        \small $\mathcal{O}(h^4)$
      }
      \end{loglogaxis}
      \node[anchor=north west] at (0.1,4.6) {(a) TraceFEM};
  \end{tikzpicture}}
\end{subfigure}
\begin{subfigure}[b]{0.49\textwidth}
  \scalebox{0.80}{\begin{tikzpicture}
    \begin{loglogaxis}[ xlabel={Mesh size $h$},domain={0.015625:0.29},
      width=1.3\linewidth, height=\linewidth,
      legend style={at={(0.99,0.015)},anchor=south east,legend columns=3,draw=none, fill=none},
      legend cell align=left, cycle list name=mark list ]
      \addlegendimage{empty legend}
      \addplot[mark=triangle] table[x=h, y=L2(P-Ph)____] {SFEM_TH_k2_apple.dat};
      \addplot[mark=square] table[x=h, y=L2(P-Ph)____] {SFEM_SF_k2_apple.dat};

      \addlegendimage{empty legend}
      \addplot[mark=triangle*] table[x=h, y=L2(P-Ph)____] {SFEM_TH_k3_apple.dat};
      \addplot[mark=square*] table[x=h, y=L2(P-Ph)____] {SFEM_SF_k3_apple.dat};

      \addplot[loosely dotted,line width=0.75pt] {1*x^2};
      \addplot[dotted,line width=0.75pt] {1*x^3};
      \addplot[dashed,line width=0.75pt] {1*x^4};
      \legend{
        ~,\small TH, \small SF ($k=2$),
        ~,\small TH, \small SF ($k=3$),
        \small $\mathcal{O}(h^2)$,
        \small $\mathcal{O}(h^3)$,
        \small $\mathcal{O}(h^4)$
      }
      \end{loglogaxis}
      \node[anchor=north west] at (0.1,4.6) {(b) SFEM};
  \end{tikzpicture}}
\end{subfigure}%
\caption{$\| p_h-p \|_{L^2(\Gamma_h)}$ errors for $k=2$ and $k=3$ for the different methods.}\label{fig:biconcave-shape_P_L2}
\end{figure}

\begin{figure}[ht]
  \centering
\begin{subfigure}[b]{0.49\textwidth}
  \scalebox{0.80}{\begin{tikzpicture}
    \def\vara{10}
    \def\varb{200}
    \begin{loglogaxis}[ xlabel={Mesh size $h$}, domain={0.015625:0.29},
      width=1.3\linewidth, height=\linewidth,
      legend style={at={(0.99,0.015)},anchor=south east,legend columns=3,draw=none, fill=none},
      legend cell align=left, cycle list name=mark list ]
      \addlegendimage{empty legend}
      \addplot[mark=triangle] table[x=h, y=H1(P-Ph)____] {TH_k2_apple.dat};
      \addplot[mark=square] table[x=h, y=H1(P-Ph)____] {SF_k2_apple.dat};
      \addlegendimage{empty legend}
      \addplot[mark=triangle*] table[x=h, y=H1(P-Ph)____] {TH_k3_apple.dat};
      \addplot[mark=square*] table[x=h, y=H1(P-Ph)____] {SF_k3_apple.dat};
      \addplot[loosely dotted,line width=0.75pt] {200*x^1};
      \addplot[dotted,line width=0.75pt] {400*x^2};
      \addplot[dashed,line width=0.75pt] {600*x^3};
      \legend{
        ~,\small TH,
        \small SF ($k=2$),
        ~,\small TH,
        \small SF ($k=3$),
        \small $\mathcal{O}(h^1)$,
        \small $\mathcal{O}(h^2)$,
        \small $\mathcal{O}(h^3)$
      }
      \end{loglogaxis}
      \node[anchor=north west] at (0.1,4.6) {(a) TraceFEM};
  \end{tikzpicture}}
\end{subfigure}
\begin{subfigure}[b]{0.49\textwidth}
  \scalebox{0.8}{\begin{tikzpicture}
    \begin{loglogaxis}[xlabel={Mesh size $h$}, domain={0.015625:0.29},
      width=1.3\linewidth, height=\linewidth,
      legend style={at={(0.99,0.015)},anchor=south east,legend columns=3,draw=none, fill=none},
      legend cell align=left, cycle list name=mark list ]
      \addlegendimage{empty legend}
      \addplot[mark=triangle] table[x=h, y=H1(P-Ph)____] {SFEM_TH_k2_apple.dat};
      \addplot[mark=square] table[x=h, y=H1(P-Ph)____] {SFEM_SF_k2_apple.dat};

      \addlegendimage{empty legend}
      \addplot[mark=triangle*] table[x=h, y=H1(P-Ph)____] {SFEM_TH_k3_apple.dat};
      \addplot[mark=square*] table[x=h, y=H1(P-Ph)____] {SFEM_SF_k3_apple.dat};

      \addplot[loosely dotted,line width=0.75pt] {5*x^1};
      \addplot[dotted,line width=0.75pt] {5*x^2};
      \addplot[dashed,line width=0.75pt] {5*x^3};
      \legend{
        ~,\small TH,
        \small SF ($k=2$),
        ~,\small TH,
        \small SF ($k=3$),
        \small $\mathcal{O}(h^1)$,
        \small $\mathcal{O}(h^2)$,
        \small $\mathcal{O}(h^3)$
      }
      \end{loglogaxis}
      \node[anchor=north west] at (0.1,4.6) {(b) SFEM};
  \end{tikzpicture}}
\end{subfigure}%
\caption{$\|\nablaGh \left( p_h-p \right) \|_{L^2(\Gamma_h)}$ errors for $k=2$ and $k=3$ for the different methods.}\label{fig:biconcave-shape_P_H1}
\end{figure}

 In Figures \ref{fig:biconcave-shape_P_L2} and \ref{fig:biconcave-shape_P_H1} the $L^2$-norm and $H^1$-norm of the pressure errors are shown. It turns out that (for both methods) the rate of convergence is less regular than for the velocity error. For both methods the $L^2$-error in the Taylor-Hood formulation converges faster than the optimal (asymptotic) convergence order $k$.  For the stream function formulations, we observe the expected optimal order $k+1$ for the $L^2$-error. The $H^1$-error in pressure for the Taylor-Hood formulation has the optimal convergence order 1 for $k=2$, but shows higher than second order convergence for $k=3$.   For the stream function formulation, in both TraceFEM and SFEM we (essentially) observe the optimal order $k$ for the $H^1$-error in pressure.

\begin{figure}[ht]
  \centering
\begin{subfigure}[b]{0.49\textwidth}
  \scalebox{0.80}{\begin{tikzpicture}
    \begin{loglogaxis}[ xlabel={Mesh size $h$},domain={0.015625:0.29},
      width=1.3\linewidth, height=\linewidth,
      legend style={at={(0.99,0.015)},anchor=south east,legend columns=2,draw=none, fill=none},
      legend cell align=left, cycle list name=mark list ]
       \addplot[mark=triangle] table[x=h, y=L2(uh*n_vbn)] {TH_k2_apple.dat};
      \addplot[mark=square] table[x=h, y=L2(uh*n_vbn)] {SF_k2_apple.dat};
      \addplot[mark=triangle*] table[x=h, y=L2(uh*n_vbn)] {TH_k3_apple.dat};
      \addplot[mark=square*] table[x=h, y=L2(uh*n_vbn)] {SF_k3_apple.dat};
      \addplot[dotted,line width=0.75pt] {1500*x^3};
      \addplot[dashed,line width=0.75pt] {25*x^4};
      \legend{
        \small TH,
        \small SF ($k=2$),
        \small TH,
        \small SF ($k=3$),
        \small $\mathcal{O}(h^3)$,
        \small $\mathcal{O}(h^4)$
      }
  \end{loglogaxis}
  \node[anchor=north west] at (0.1,4.6) {(a) TraceFEM};
  \end{tikzpicture}}
\end{subfigure}
\begin{subfigure}[b]{0.49\textwidth}
  \scalebox{0.80}{\begin{tikzpicture}
    \begin{loglogaxis}[ xlabel={Mesh size $h$},domain={0.015625:0.29},
      width=1.3\linewidth, height=\linewidth,
      legend style={at={(0.99,0.015)},anchor=south east,legend columns=2,draw=none, fill=none},
      legend cell align=left, cycle list name=mark list ]
       \addplot[mark=triangle] table[x=h, y=L2(uh*n_vbn)] {SFEM_TH_k2_apple.dat};
      \addplot[mark=square] table[x=h, y=L2(uh*n_vbn)] {SFEM_SF_k2_apple.dat};
      \addplot[mark=triangle*] table[x=h, y=L2(uh*n_vbn)] {SFEM_TH_k3_apple.dat};
      \addplot[mark=square*] table[x=h, y=L2(uh*n_vbn)] {SFEM_SF_k3_apple.dat};
      \addplot[dotted,line width=0.75pt] {40*x^3};
      \addplot[dashed,line width=0.75pt] {0.1*x^4};
      \legend{
        \small TH, \small SF ($k=2$),
        \small TH, \small SF ($k=3$),
        \small $\mathcal{O}(h^3)$,
        \small $\mathcal{O}(h^4)$
      }
      \end{loglogaxis}
      \node[anchor=north west] at (0.1,4.6) {(b) SFEM};
  \end{tikzpicture}}
\end{subfigure}%
\caption{$\|\bu_h \cdot \tilde \bn_h\|_{L^2(\Gamma_h)}$ errors for $k=2$ and $k=3$ for the different methods.}\label{fig:biconcave-shape_U_normal}
\end{figure}

We further determine $\|\bu_h\cdot\tilde{\bn}_h\|_{L^2(\Gamma_h)}$ to measure how well the numerical solution satisfies the tangential condition. Note that we use the improved normal vector $\tilde{\bn}_h$ that is used in the penalty term in the Taylor-Hood formulation and in the velocity reconstruction in the stream function formulation. In Figure \ref{fig:biconcave-shape_U_normal} we see an order $k+1$  convergence for both, the Taylor-Hood and stream function formulation in both discretization methods. Note that in the stream function formulation on the continuous level, due to  the relation $\bu^\ast_T= \bn \times \gradG \psi^\ast$, the tangential condition is automatically fulfilled (cf. section~\ref{sectdisc}). This explains why in   Figure \ref{fig:biconcave-shape_U_normal} the quantity $\|\bu_h\cdot\tilde{\bn}_h\|_{L^2(\Gamma_h)}$ is (much) smaller for the stream function formulation than for the Taylor-Hood formulation.

\begin{figure}[ht]
  \centering
\begin{subfigure}[b]{0.49\textwidth}
  \scalebox{0.80}{\begin{tikzpicture}
    \begin{loglogaxis}[ xlabel={Mesh size $h$},domain={0.015625:0.29},
      width=1.3\linewidth, height=\linewidth,
      legend style={at={(0.99,0.015)},anchor=south east,legend columns=2,draw=none, fill=none},
      legend cell align=left, cycle list name=mark list ]
      \addplot[mark=triangle] table[x=h, y=L2(DivG_uh)_] {TH_k2_apple.dat};
      \addplot[mark=square] table[x=h, y=L2(DivG_uh)_] {SF_k2_apple.dat};
      \addplot[mark=triangle*] table[x=h, y=L2(DivG_uh)_] {TH_k3_apple.dat};
      \addplot[mark=square*] table[x=h, y=L2(DivG_uh)_] {SF_k3_apple.dat};
      \addplot[dotted,line width=0.75pt] {750*x^2};
      \addplot[dashed,line width=0.75pt] {300*x^3};
      \legend{
        \small TH, \small SF ($k=2$),
        \small TH, \small SF ($k=3$),
        \small $\mathcal{O}(h^2)$,
        \small $\mathcal{O}(h^3)$
      }
      \end{loglogaxis}
      \node[anchor=north west] at (0.1,4.6) {(a) TraceFEM};
  \end{tikzpicture}}
\end{subfigure}
\begin{subfigure}[b]{0.49\textwidth}
  \scalebox{0.80}{\begin{tikzpicture}
    \begin{loglogaxis}[ xlabel={Mesh size $h$},domain={0.015625:0.29},
      width=1.3\linewidth, height=\linewidth,
      legend style={at={(0.99,0.015)},anchor=south east,legend columns=2,draw=none, fill=none},
      legend cell align=left, cycle list name=mark list ]
      \addplot[mark=triangle] table[x=h, y=L2(DivG_uh)_] {SFEM_TH_k2_apple.dat};
      \addplot[mark=square] table[x=h, y=L2(DivG_uh)_] {SFEM_SF_k2_apple.dat};
      \addplot[mark=triangle*] table[x=h, y=L2(DivG_uh)_] {SFEM_TH_k3_apple.dat};
      \addplot[mark=square*] table[x=h, y=L2(DivG_uh)_] {SFEM_SF_k3_apple.dat};
      \addplot[dotted,line width=0.75pt] {15*x^2};
      \addplot[dashed,line width=0.75pt] {3*x^3};
      \legend{
        \small TH, \small SF ($k=2$),
        \small TH, \small SF ($k=3$),
        \small $\mathcal{O}(h^2)$,
        \small $\mathcal{O}(h^3)$
      }
      \end{loglogaxis}
      \node[anchor=north west] at (0.1,4.6) {(b) SFEM};
  \end{tikzpicture}}
\end{subfigure}%
\caption{$\| \DivGh \left( \bu_h \right) \|_{L^2(\Gamma_h)}$ errors for $k=2$ and $k=3$ for the different methods.}\label{fig:biconcave-shape_U_divergence}
\end{figure}

Since we are interested in solenoidal vector fields $\bu$, the error in the divergence, $\DivGh \left( \bu_h \right)$, is measured and plotted in Figure \ref{fig:biconcave-shape_U_divergence}. In both, the discrete solution of the Taylor-Hood formulation and the reconstructed velocity vector in the stream function formulation, the convergence order $k$ can clearly be observed in these results.

We conclude that with the parameter settings (penalty parameter, stabilization parameter) defined above \emph{the two methods, TraceFEM and SFEM, have (essentially) the same rate of convergence} for all considered quantities: velocity and pressure errors, tangentiality measure and discrete divergence. For the Taylor-Hood formulation the optimal velocity error convergence orders are $k$ and $k+1$ for the $H^1$- and $L^2$-norm, respectively, and the optimal orders for the pressure error are $k-1$ and $k$ for the $H^1$- and $L^2$-norm, respectively. For the stream function formulation the optimal orders are $k$ and $k+1$ for the $H^1$- and $L^2$-norm of the velocity error, and $k$ and $k+1$ for the $H^1$- and $L^2$-norm of the pressure error. The numerical results show that \emph{these optimal orders are attained} and in certain cases, e.g., the pressure $L^2$-norm error in the Taylor-Hood formulation, we even observe a higher rate than the  optimal one. For the tangentiality measure $\|\bu_h \cdot \tilde \bn_h\|_{L^2(\Gamma_h)}$ we obtain the convergence order $k+1$ in both the Taylor-Hood and the stream function formulation.  For the discrete divergence $\|\DivGh \left( \bu_h \right)\|_{L^2(\Gamma_h)}$ we obtain the convergence order $k$ in both the Taylor-Hood and the stream function formulation.
We should remark that optimal error convergence orders for SFEM are not known analytically and are here considered as expected orders based on the results for a vector surface Laplacian \cite{hansbo2019analysis}. Other open issues, which are postponed to future investigations, are error norms for the tangential velocity.

\subsection{Discussion of results}
We summarize and discuss a few aspects of the methods treated above.\\
\emph{Surface approximation.} In both the TraceFEM and SFEM one needs an initial piecewise planar approximation $\Gamma^{\rm lin}$ of $\Gamma$ that is sufficiently accurate, in the sense that the closest point projection $\pi:\Gamma^{\rm lin} \to \Gamma$ should be a bijection and ${\rm dist}(\Gamma^{\rm lin}, \Gamma) \sim h^2$. For the TraceFEM it is easy to construct such a $\Gamma^{\rm lin}$, using linear finite interpolation (or approximation) of the level set function on a volume triangulation. For SFEM, opposite to TraceFEM, the approximation $\Gamma^{\rm lin}$ has to be shape-regular and in general the construction of $\Gamma^{\rm lin}$ is more difficult, in particular for surfaces with strongly varying curvatures. Given $\Gamma^{\rm lin}$, a higher order surface approximation is obtained by a suitable parametric approach. In TraceFEM this is based on a transformation (deformation) of the local volume triangulation $\cT_h^\Gamma$ , whereas in SFEM a transformation of the surface approximation $\Gamma^{\rm lin}$ is used. \\
\emph{Stream function formulation or formulation in $(\bu,p)$-variables.} The stream function formulation can only be used if $\Gamma$ is simply connected. Based on the complexity results and the error plots presented above, we conclude that (for this test case)  the methods based on the stream function formulation are for $k \leq 3$ slightly more efficient than the ones based on the $(\bu,p)$-variables. There is, however, not a decisive difference in efficiency. In the stream function formulation the tangential contraint is automatically satisfied, whereas in  the $(\bu,p)$ formulations we need a penalty approach. \\
\emph{Complexity}. If one defines complexity in terms of DOFs, cf. Section~\ref{sectC}, we obtain  complexity estimates $\sim n k^3$  and $\sim  n k^2$ (precise relations in \eqref{nr1}--\eqref{nr3}) for the TraceFEM and SFEM, respectively. The difference in the exponents is caused by the fact that the TraceFEM uses finite element spaces on a strip of 3D tetrahedral elements whereas the SFEM is based on spaces on 2D triangular elements.  \\
\emph{Rate of convergence and efficiency.} In the numerical experiments we obtain for both TraceFEM and SFEM  in the $(\bu,p)$ formulation and the stream function formulation optimal  rates of convergence in $L^2$- and $H^1$ norms. Looking at the size of the errors we see that for comparable DOFs the SFEM typically has an error that is 10-50 times smaller than the corresponding error in the TraceFEM. Hence, if one uses error/DOFs as measure of efficiency then the SFEM is significantly more efficient than TraceFEM. \\
\emph{Parameter tuning.} The parameters are summarized in Table~\ref{tabpar}. Note that for SFEM in
stream function formulation there are no parameters. The use of additional scaling constants, e.g., $c_p, c_\eta$ in $\rho_p= c_p h$, $\eta=c_\eta h^{-2}$, did note significantly change the results. \\
\emph{Applications to other problem classes.} The methods  studied in this paper can also be used in other (more complex) surface PDE problems. For example, the surface Stokes problem may be coupled to a bulk flow problem. In such a setting the TraceFEM has the advantage that for the bulk and surface problems one can use the same finite element spaces. Another extension concerns (Navier-)Stokes equations on evolving surfaces. If the surface evolution is smooth and without strong deformations we expect, based on the results presented in this paper, that the evolving SFEM, cf. \cite{DE2007Finite} for scalar-valued problem, is more efficient than a time dependent variant of TraceFEM.  If, however, the surface geometry has strong deformations, or even topological singularities occur, the TraceFEM may be more attractive.

\subsection{Benchmark problem}
\begin{figure}[ht!]
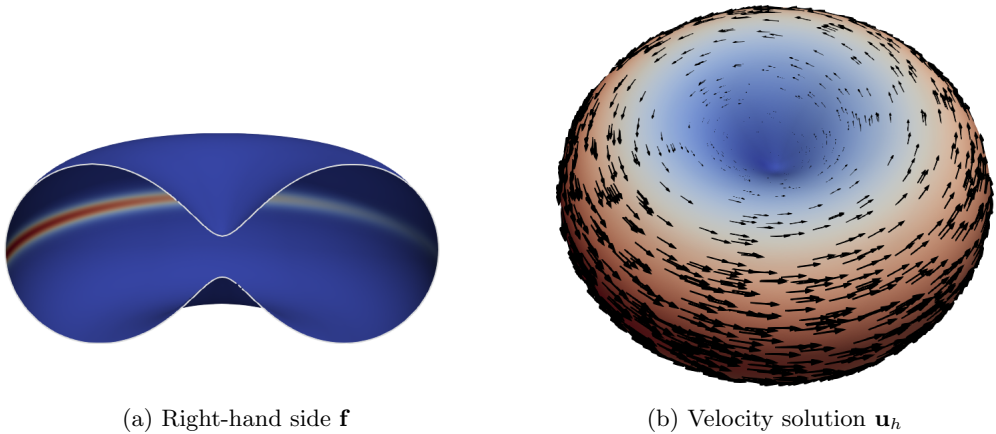

  \centering
  \begin{subfigure}[b]{0.45\textwidth}
  \hspace*{-0.5cm}
      \centering
      \includegraphics[width=\textwidth]{figure_9_0.png}
      \caption{Right-hand side $\bbf$}
  \end{subfigure}
  \hfill
  \begin{subfigure}[b]{0.45\textwidth}
      \centering
      \includegraphics[width=\textwidth]{figure_9_1.png}
      \caption{Velocity solution $\bu_h$}
  \end{subfigure}
  \caption{Geometry ``biconcave shape'' for $d=0.96$ in benchmark problem. The right-hand side is  non-zero only on  a thin strip located close to the $y$-$z$-plane.  Glyphs represents the force and velocity magnitude by  coloring from blue (small) to red (large) (Color coding available online).}\label{fig:benchmark_rhs_u}
\end{figure}

We consider the Stokes problem on the biconcave shape as in \eqref{eq:biconcave-shape_levelset}, but with a prescribed
right-hand side function $\bbf$. This driving force is constructed in such a way that a rotating flow around the $x$-axis emerges (cf. Fig.~\ref{fig:biconcave-shape}(b) for axes-directions).   Choosing $\bbf$ asymmetric places the two emerging vortices away from the center position at the geometric bumps on the $x$-axis, cf. Fig.~\ref{fig:benchmark_rhs_u}.  The position of the two vortices depends on the size of the curvature at the bumps, which is controlled by the geometry parameter $d$ in \eqref{eq:biconcave-shape_levelset}. This geometric effect has already been discussed in \cite{ReVo2015,ReVo2018} for the surface Navier-Stokes equations and in \cite{TuViNe2010} for surface super fluids. The latter case allows for analytic expressions for the interaction of vortices with the surface geometry. We here only use the geometry effect to formulate a benchmark problem for the surface Stokes equations with the center of the vortices as quantity of interest.

For the construction of the right-hand side function, we consider the rotational tangential field $\bbf_0 \colonequals \bn\times (1,0,0)^T$. This tangential vector field is first restricted to an outer ring of the surface and then accelerated depending on the rotational angle, as follows:
\begin{equation*}
  \bbf(x) \colonequals \chi_\varepsilon(\xx)\,\tfrac{1}{2}\big(1 + \sin(\alpha(\xx))\big)\cdot\bbf_0(\xx)
\end{equation*}
with $\chi_\varepsilon(\xx) \colonequals \delta_\varepsilon(x)\,\delta_\varepsilon(\sqrt{y^2 + z^2} - R)$, $\alpha(\xx) \colonequals \arctan\!2(y,z)$,
\begin{equation*}
  \delta_\varepsilon(r) \colonequals 36\varphi_\varepsilon(r)^2\big(1-\varphi_\varepsilon(r)\big)^2\,,\text{ and }
  \varphi_\varepsilon(r) \colonequals \tfrac{1}{2}\big(1 - \tanh(3\,r/\varepsilon)\big)\,.
\end{equation*}

This construction follows general ideas of \cite{LiLoRaVo2009}. The right-hand side $\bbf$ and the velocity solution $\bu_h$ are illustrated for $d=0.96$ in Fig.~\ref{fig:benchmark_rhs_u}. For the parameters we have used $R=1.1$ for the outer ring radius and $\varepsilon=0.2$ for the restriction thickness. The angle-dependent scaling factor is chosen such that it has a maximum on one side and the minimum zero on the opposite site of the shape.

For this given right-hand side, the Stokes equations are solved using the four approaches discussed before, i.e., the TraceFEM and SFEM discretization of the Taylor-Hood and stream function formulation, respectively. In the Taylor-Hood formulation, for the resulting velocity field $\bu_h$, inside a vortex we determine the
location where the size of the velocity is minimal and take this as approximation of the vortex position. For the  stream function approximation $\psi_h$, the  location of a local maximum yields the numerical approximation of the vortex position. Below this numerical approximation of the vortex position is denoted by $\bx_v$.

\begin{figure}[ht!]
  \centering
  \begin{subfigure}[b]{0.24\textwidth}
      \centering
      \includegraphics[width=\textwidth]{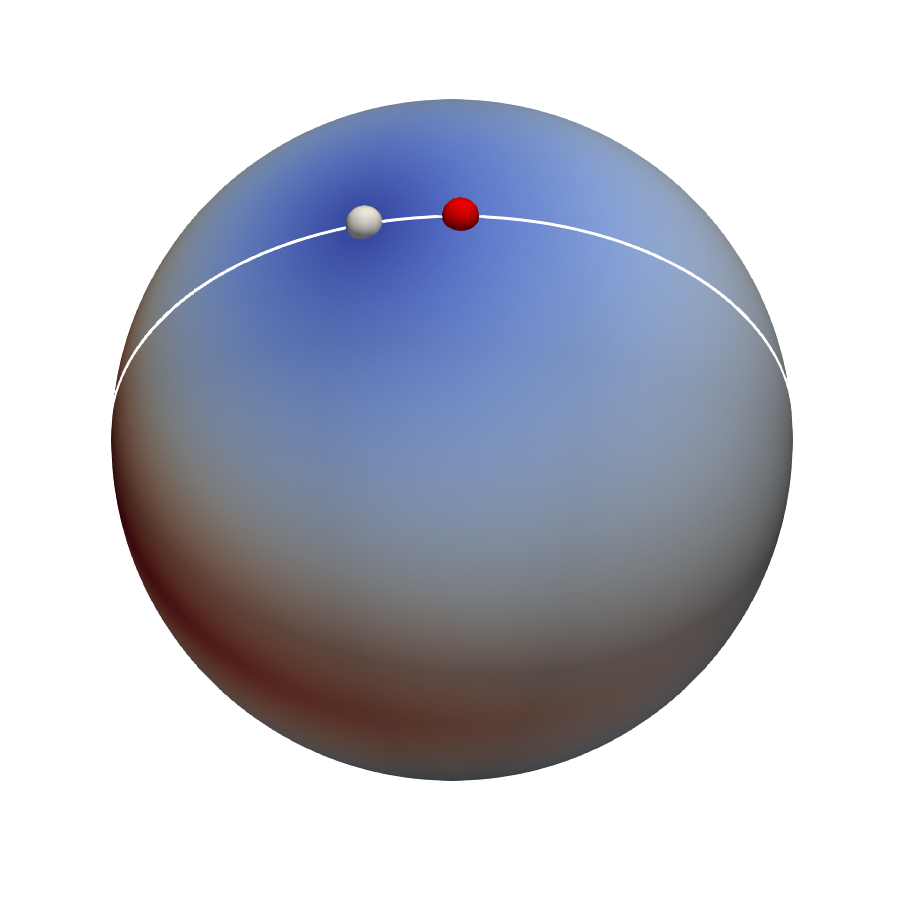}
      \caption{$d=0.0$}
  \end{subfigure}
  \hfill
  \begin{subfigure}[b]{0.24\textwidth}
      \centering
      \includegraphics[width=\textwidth]{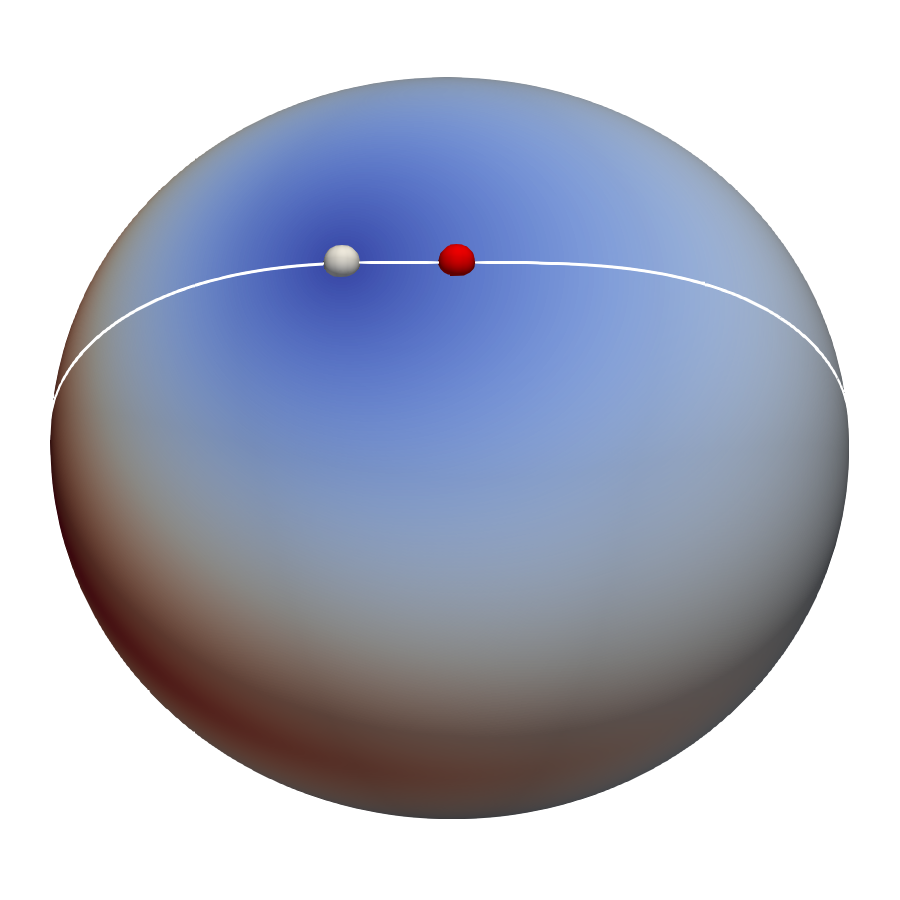}
      \caption{$d=d_0\approx 0.572$}
  \end{subfigure}
  \hfill
  \begin{subfigure}[b]{0.24\textwidth}
      \centering
      \includegraphics[width=\textwidth]{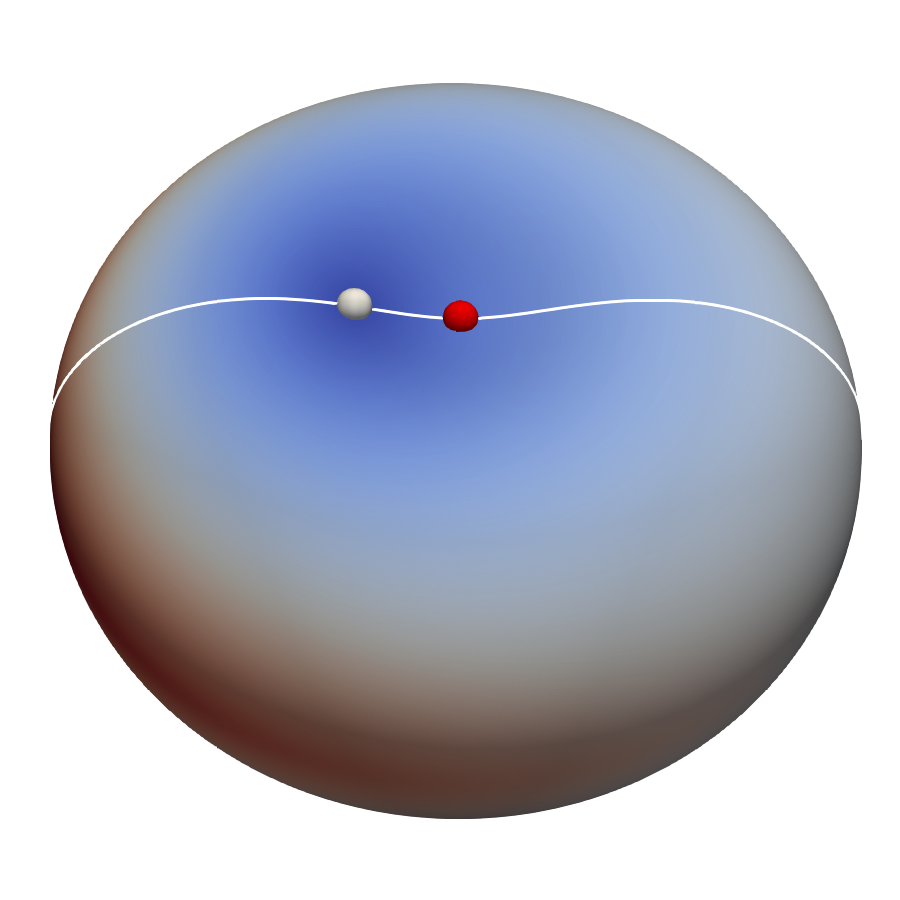}
      \caption{$d=0.8$}
  \end{subfigure}
  \hfill
  \begin{subfigure}[b]{0.24\textwidth}
      \centering
      \includegraphics[width=\textwidth]{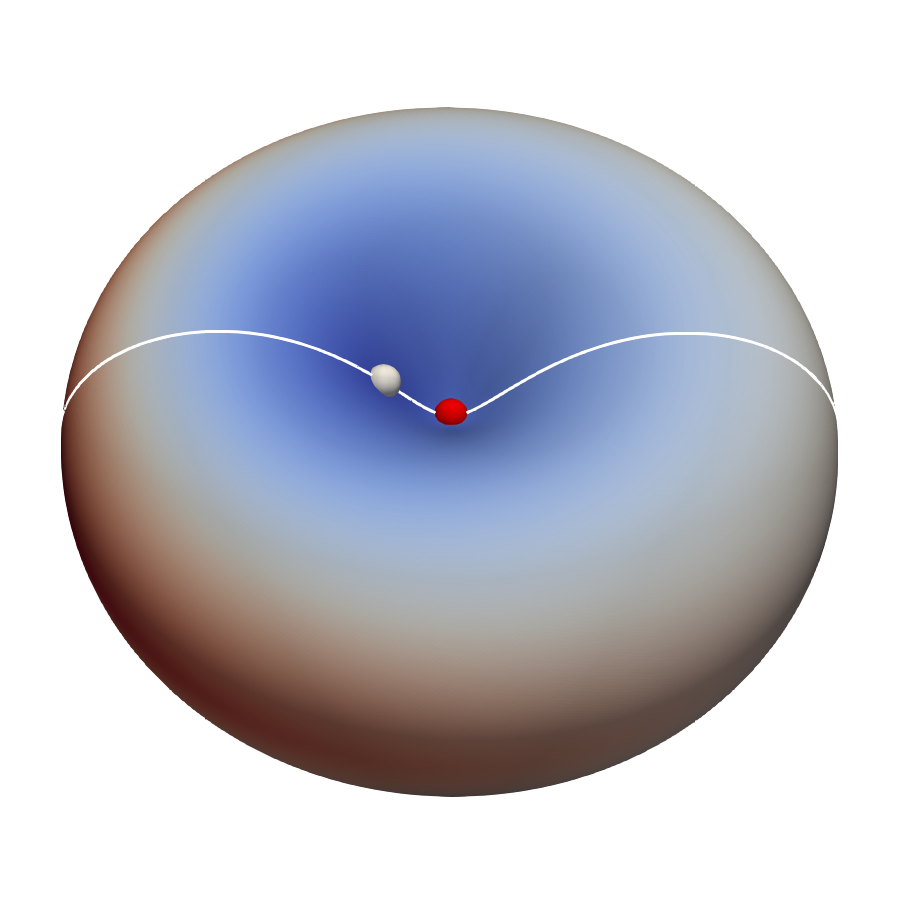}
      \caption{$d=0.96$}
  \end{subfigure}
  \caption{Magnitude of the computed velocity field $\bu_h$ for different shape parameters $d$. In white  numerical vortex location and in red the center of the geometry. The slice highlights the direction of the asymmetry of the force field and the shape of the geometry in this direction (Color coding available online).}\label{fig:biconcave-shape_u_vertices}
\end{figure}

In Figure \ref{fig:biconcave-shape_u_vertices}, the velocity magnitude and the location of the vortex on the upper side ($x > 0$) is visualized. The Euclidean distance of the vortex location $\bx_v$ to the center on the surface $\bx_c=(\sqrt{ c^{4/3} - d^2},0,0)^T$ and the dependence of this distance on the geometric parameter $d$ are the quantities of interest. We study four different cases, namely $d=0$, $d=d_0\colonequals\sqrt{\frac{3}{8} c^{8/3}}\approx 0.572$, $d=0.8$, and $d=0.96$. The value $d_0$ corresponds to a geometry with zero mean and Gaussian curvature in the geometry center $\bx_c$. The  values for $d$ with corresponding mean and Gaussian curvature in the geometry center $\bx_c$ are given in Table \ref{tab:Gaussian_curvature_x_c}.
\begin{table}[ht!]
  \begin{center}
  \begin{tabular}{c|c|c|c|c}
    \toprule
    & $d=0$ & $d=d_0$  & $d=0.8$ & $d=0.96$ \\
    \hline
    $K$     & $1.07$ & $0$ & $3.12$ & $268.76$ \\
    \hline
    $\tr(\bH)$    & $2.07$ & $0$ & $-3.53$ & $-32.79$ \\
    \bottomrule
  \end{tabular}
  \end{center}
  \caption{Mean curvature $\tr(\bH)$ and Gaussian curvature $K$ evaluated at the geometry center for various geometry parameters $d$.}\label{tab:Gaussian_curvature_x_c}
\end{table}

We performed numerical experiment in the same setting as explained in Section~\ref{sec:numerical_experiments}. From the results presented in that section we see that on the finest level 5 and with $k=3$ the most accurate results in all cases (except the $H^1$-norm velocity error in TraceFEM) are obtained using the stream function formulation. Compare Table~\ref{tab:comparison_data-2} for the actual complexity data of the grids for the geometry $d=0.96$ used in SFEM and TraceFEM. For the other geometry parameters the number of elements is in the same order.

\begin{table}[ht]
  \begin{center}
  \begin{tabular}{r|l|l|l|l}
    \toprule
               & grid size   & grid elements     &  nnz & DOFs \\
    \hline
    SFEM       & $h\;=0.009984$      & $528384$ & $5.59\cdot 10^8$ & $8454148$ \\
    TraceFEM   & $h_\Gamma=0.009744$ & $411354$ & $7.52\cdot 10^8$ & $11065384$ \\
    \bottomrule
  \end{tabular}
  \end{center}
  \caption{Complexity of SFEM and TraceFEM grids for the geometry parameter $d=0.96$ and grid refinement level $5$: (Surface) grid size, number of elements in the computational grid,  number of non-zeros (nnz) and number of degrees of freedom (DOFs) for $k=3$ in the stream function formulation. Compare also Table~\ref{tab:comparison_data}}\label{tab:comparison_data-2}
\end{table}

In Table \ref{tab:vortex_distance_ref} we present the distance results $\|\bx_v - \bx_c\|_2$ both for the TraceFEM and SFEM discretization of the stream function formulation for $k = 3$ and refinement level 5.

\begin{table}[ht!]
  \begin{center}
  \begin{tabular}{r|l|l|l|l}
    \toprule
    & $d=0$ & $d=d_0$  & $d=0.8$ & $d=0.96$ \\
    \hline
    SFEM     & $0.255577$ & $0.308290$ & $0.295497$ & $0.245279$ \\
    TraceFEM & $0.254524$ & $0.309088$ & $0.295475$ & $0.244346$ \\
    \bottomrule
  \end{tabular}
  \end{center}
  \caption{Reference solution of the distance of the vortex center to the geometry center for various geometry parameters $d$.}\label{tab:vortex_distance_ref}
\end{table}

Note that we extensively tested the stream function formulation and compared it with the  Taylor-Hood formulation (previous section) and that we use two different methods (TraceFEM and SFEM) that are implemented in two different software codes. Based on this we claim that the first 3-4 digits of the distance values in Table~\ref{tab:vortex_distance_ref} are correct. These results can be used as \emph{benchmark values} for the development and testing of other codes used for the numerical simulation of surface Stokes equations.

\section*{Acknowledgements}
The authors wish to thank the German Research Foundation (DFG) for financial support within the Research Unit ``Vector- and Tensor-Valued Surface PDEs'' (FOR 3013) with project no. RE 1461/11-1 and VO 899/22-1. We further acknowledge computing resources provided by ZIH at TU Dresden and within project PFAMDIS at FZ J{\"u}lich.

\bibliographystyle{siam}
\bibliography{bibliography}
\end{document}